\documentclass[11 pt, psamsfonts]{amsart}
\usepackage{amssymb}
\usepackage{array,calc,url,tabularx, epic,tikz}
\usetikzlibrary{arrows,positioning,decorations.pathmorphing,
  decorations.markings}

\def\ignore#1{\relax}

\def\g{\mathfrak g}
\def\q{\mathfrak q}
\def\m{\mathfrak m}

\def\Z{{\mathbb Z}}

\def\C{{\mathbb C}}

\def\Ga{\Gamma}

\def\la{\lambda}
\def\La{\Lambda}

\def\H{\mathcal H}

\def\A{\mathcal A}
\def\Ca{\mathcal C}

\def\H{\mathcal H}
\def\P{\mathcal P}

\def\A{\mathcal A}

\def\K{\overline{K}}
\def\ni{\noindent}

\def\ignore#1{\relax}

\def\1{{\bf 1}}

\def\ep{\epsilon}

\def\End{{\rm End}}

\def\ve{\varepsilon}

\def\p{{\frak p}}

\def\ni{\noindent}

\calclayout

\renewenvironment{proof}[1][\proofname]{\par
  \normalfont
  \topsep6\p@\@plus6\p@ \trivlist
  \item[\hskip\labelsep\bfseries
    #1\@addpunct{.}]\ignorespaces
}{%
  \qed\endtrivlist
}

\def\th@plain{%
  \let\thmhead\thmhead@plain \let\swappedhead\swappedhead@plain
  \thm@preskip.5\baselineskip\@plus.2\baselineskip
                                    \@minus.2\baselineskip
  \thm@postskip\thm@preskip
  \itshape
\renewcommand{\labelenumi}{{(\alph{enumi})\quad}}
                        \renewcommand{\labelenumii}{{(\roman{enumii})\ }}
}
\def\th@definition{%
  \let\thmhead\thmhead@plain \let\swappedhead\swappedhead@plain
  \thm@preskip.5\baselineskip\@plus.2\baselineskip
                                    \@minus.2\baselineskip
  \thm@postskip\thm@preskip
  \upshape
}
\def\th@remark{%
  \thm@headfont{\itshape}
  \let\thmhead\thmhead@plain \let\swappedhead\swappedhead@plain
  \thm@preskip.5\baselineskip\@plus.2\baselineskip
                                    \@minus.2\baselineskip
  \thm@postskip\thm@preskip
  \upshape
}

{\theoremstyle{plain}
\newtheorem{theorem}{Theorem}[section]
}

{\theoremstyle{plain}
\newtheorem{proposition}[theorem]{Proposition}
}

{\theoremstyle{plain}
\newtheorem{corollary}[theorem]{Corollary}
}

{\theoremstyle{plain}
\newtheorem{lemma}[theorem]{Lemma}
}

{\theoremstyle{plain}

}

{\theoremstyle{definition}
\newtheorem{definition}[theorem]{Definition}
}

{\theoremstyle{definition}
\newtheorem{example}[theorem]{Example}
}

{\theoremstyle{remark}
\newtheorem{remark}[theorem]{Remark}
}

{\theoremstyle{remark}

}

\numberwithin{equation}{section}

\renewcommand{\labelenumi}{{ \theenumi.}}
\renewcommand{\labelenumii}{{(\alph{enumii})}}

\def\ni{\noindent}

\def\la{\lambda}
\def\al{\alpha}

\def\ga{\gamma}

\def\choose #1 #2{\begin{pmatrix}#1\\#2\end{pmatrix}}

\def\K{{\mathcal K}}
\def\A{{\mathcal A}}
\def\H{{\mathcal H}}
\def\BM{{\mathcal B}{\mathcal M}}

\begin{document}

\title[Simple $B_4$ representations]
{Simple $B_4$ representations associated to cyclotomic Hecke algebras}

\author{Lilit Martirosyan }
\address{L.M. Department of Mathematics and Statistics\\ University of North Carolina\\ Wilmington \\North Carolina}
\email{martirosyanl@uncw.edu}

\author{Hans Wenzl}
\address{H.W. Department of Mathematics\\ University of California\\ San Diego,
California}

\email{hwenzl@ucsd.edu}

\begin{abstract} 
We determine the structure of the cyclotomic Hecke algebra corresponding to the complex reflection group $G_{25}$ also when it is not semisimple, as long as the generators are diagonalizable. In particular, we classify all simple representations of the braid group $B_4$ for which the generators are diagonalizable and satisfy a cubic polynomial. This will be used in the classification of braided tensor categories of type $G_2$. 
\end{abstract}
\maketitle

It is a classical result by Coxeter that the relation $\sigma_i^p=1$ defines a finite quotient of the braid group $B_n$
for $p>2$ and $n>2$  if and only if $p=3$ and $n\leq 5$ or $n=3$ and $p\leq 5$. The corresponding quotient
groups are special cases of complex reflection groups. Hecke algebra deformations have been defined for these
groups which, for generic parameters, have the same representation theory as the group ring. The main work in this paper 
consists of studying the representations for the Hecke algebra $\K_4$ for $n=4$ and $p=3$. It can be defined as a quotient of the group algebra of $B_4$
with the additional relation $\prod_{i=1}^3 (\sigma_j-\la_i1)$ and corresponds to the
complex reflection group with the label $G_{25}$. While this seems to be a worthwhile problem in its own
right, our research was motivated by our work on classifying tensor categories with the tensor product rules
of the simple Lie algebra of type $G_2$. We will only need representations for which the image of the braid
generators are diagonalizable. So our results will mostly be about this case, which is also more accessible to the methods we use.
Our main results in this paper are:
\vskip .2cm
\ni {\bf Theorem A} (see Theorem \ref{semisimple:thm}). The algebra $\K_4$ is not semisimple if and only if the eigenvalues of $\sigma_1$ are zeros of one of the following polynomials: $\la_i-\la_j$,  $\la_i+\la_j$, $\la_i+\theta\la_j$, 
$\la_i\pm\sqrt{-1}\la_j$, $\la_i^2+\la_j\la_k$, $\la_i^2-\theta\la_j\la_k$ or $\la_i^3-\la_j^2\la_k$, where $\{ i,j,k\}=\{ 1,2,3\}$ and $\theta$ is a primitive third root of unity.
\vskip .2cm
This result would also follow from the calculation of Schur elements of a conjectured symmetrizing trace on $\K_4$ in \cite{Malledegree}, see the discussion in Section \ref{Schursec}.
There exists a complete set of mutually non-isomorphic $\K_4$ modules $V_\ga$ each of which is simple for a Zariski open set of eigenvalues; we call these regular modules. The canonical central element $\Delta_4^2$ acts via a scalar $\la_\ga$ on $V_\ga$ which is a monomial in the eigenvalues (possibly multiplied with a primitive third root of unity, which we ignore for the moment). Let $\p$ be the prime ideal generated by one of the polynomials $\neq \la_i-\la_j$ listed in Theorem A. We say $\ga\sim_\p\ga'$ if $\la_\ga-\la_{\ga'}\in \p$. 
\vskip .2cm
\ni {\bf Theorem B} (see Theorem \ref{oneideal})  Let $R$ be a domain in which $\la_i\neq\la_j$ for $i\neq j$.
Then the blocks of $\K_4$, defined over the quotient field $Q(R/\p)$, are given by the equivalence classes with respect to $\sim_\p$. Moreover, if we order the labels in such an equivalence class $\ga_1>\ga_2>\ ...\ \ga_r$
via alphabetical order of the monomials $\la_\ga$, we obtain an exact sequence
$$0\ \to\  V_{\ga_1}\ \to\  V_{\ga_2}\ \to\  ...\  \to\  V_{\ga_r}\ \to\  0,$$
where the induced subquotient modules for each $V_{\ga_i}$ are simple.
\vskip .2cm
For more details, see the table before Theorem \ref{oneideal}.
The simple modules in Theorem B may not stay simple if the eigenvalues satisfy more than one polynomial in Theorem A. We determine all possible simple subquotient modules also in these cases, see Table 4.
It turns out that all of these modules are isomorphic to either regular modules or modules appearing in Theorem B. Hence we  obtain 

\vskip .2cm
{\bf Theorem C} (see Theorem \ref{nomoredecomposition}) Table 1 and Table 3 give a complete classification of all simple modules of $\K_4$, subject to the condition that the braid generators are diagonalizable. 
\vskip .2cm
Unfortunately, even though we can give a fairly uniform description of our results, our proofs involve a lot of calculations and case-by-case studies. In order to study representations which are not already known from the study of Hecke algebras and BMW algebras, we use the explicit representations of $\K_4$ obtained in \cite{MM}. To determine their structure when they are not simple, we use different bases which we call weight and path bases. We then calculate diagonal matrix entries of eigenprojections of braid generators  with respect
to these bases, using a generalization of the quantum Jucys-Murphy elements. The main technical tool is Theorem \ref{AB2theorem}, which first appeared in \cite{Wexc}. We can then read off these diagonal entries
the existence and structure of subquotient modules for the representations in \cite{MM}; this is outlined at the beginning of Section \ref{structure:sec}. Here is our paper in more detail.

We review basic facts about the cubic Hecke algebras $\K_3$ and $\K_4$, regular modules and BMW algebras in Section 1.
We also introduce path and weight bases. In Section 2 we find restrictions on the eigenvalues for which $\K_4$ is not semisimple, and show that path bases exist in all such cases. Section 3 is devoted to calculate the previously mentioned matrix entries for path bases. This is then used in Section 4 to determine all blocks and all possible simple representations of $\K_4$, assuming we have no repeated eigenvalues. We discuss connections of our results to other results and possible applications in the last section.

$Acknowledgement:$ We would  like  to thank Ivan Marin and Gunter Malle for useful information. H.W. would like to thank
UNC Wilmington and Lilit Martirosyan for supporting two visits. L.M. gratefully acknowledges a Reassignment Award from the College of Science and Engineering at UNC Wilmington. Both authors would like to thank MFO Oberwolfach for supporting us
as research fellows for two weeks  in October 2024.

\section{Cyclotomic Hecke algebras for $G_4$ and $G_{25}$} 

\subsection{Definitions}\label{definitionsection} We first review a few basic facts about Artin's braid group, see e.g. \cite{KT}. The braid group $B_n$ is defined 
by generators $\sigma_i$, $1\leq i<n$ and relations $\sigma_i\sigma_j=\sigma_j\sigma_i$ for $|i-j|>1$ and $\sigma_i\sigma_{i+1}\sigma_i=\sigma_{i+1}\sigma_i\sigma_{i+1}$ for $i<n-1$. We  define elements $\Delta_k$ inductively by $\Delta_1=1$ and 
$$\Delta_{k+1}=\sigma_1\sigma_2\ ...\ \sigma_k\Delta_k=\Delta_k\sigma_k\sigma_{k-1}\ ...\ \sigma_1.$$
 Then $\Delta_n^2$ is the generator of the center of $B_n$, see e.g. \cite{KT} Section 1.3.3.
Let $\gamma_n=\Delta_n^2\Delta_{n-1}^{-2}$. Then $\ga_n=\sigma_{n-1}\ ...\ \sigma_2\sigma_1^2\sigma_2\ ...\ \sigma_{n-1}$, which can be fairly easily seen by using the usual description of braids via pictures, see e.g. \cite{KT}.
We will also need the following well-known identity
\begin{equation}\label{braidpermutation}
\Delta_4\sigma_i\Delta_4^{-1}\ =\ \sigma_{4-i},\quad 1\leq i<4.
\end{equation}
It was shown by Coxeter \cite{Cox}
that $B_n$ with the additional relation $\sigma_i^3=1$ is a finite group only if $n\leq 5$. These quotients correspond to the complex
reflection groups denoted by $G_4$, $G_{25}$ and $G_{32}$ in \cite{ST}. This was generalized to the quotient $\K_n$  of the group algebra $\C B_n$ of the braid group $B_n$
defined by the relations $(\sigma_i-\la_1)(\sigma_i-\la_2)(\sigma_i-\la_3)=0$, $1\leq i<n$ as follows.

\begin{theorem} \cite{BMR}, \cite{MHecke} The algebras $\K_n$ have dimensions 3, 24, 648 and 155,520 for $n=2,3,4,5$.  In particular, they have bases with respect to which the standard generators
act via matrices whose entries are Laurent polynomials in the eigenvalues $\la_i$ over a finite extension of $\Z$. Their $\C$ span is isomorphic to the group ring of the corresponding
complex reflection group for a Zariski open set of the eigenvalues.

The algebras $\K_n$ are
infinite-dimensional for $n>5$.
\end{theorem}
We define a certain class of modules, called regular modules,  which play the role of Specht modules (for Hecke algebras of type $A$) or cellular modules (for $BMW$ algebras, see Section \ref{BMW:sec}) for our algebras $\K_n$.
\begin{definition}\label{regularrep} We call a $\K_n$ module (for $n\leq 5$) a $regular$ module, if 

(a) it has a basis with respect to which  the braid generators  are represented by matrices whose
entries are in a finite extension of $\Z[\la_1,\la_2,\la_3]$,

(b) the representation is irreducible for a Zariski open subset of the eigenvalues $(\la_i)$.
\end{definition}
\begin{remark}\label{pathunique}
It is easy to see that two regular modules which are isomorphic for a Zariski dense subset of the eigenvalues may not be isomorphic when they are no longer simple. 
Indeed, as any braid group $B_n$ is isomorphic to its opposite group,
we obtain from any representation $\rho: B_n\to \End(V)$  a representation $\rho^t$ defined by $\rho^t(\beta)=\rho(\beta)^t$,
the transpose matrix of $\rho(\beta)$. So any quotient module of $V$ corresponds to a submodule of $V^t$ (which is $V$ with the $\rho^t$ action of $B_n$).
 If $V$ is not semisimple, this entails that $V^t$ is usually not isomorphic to $V$. More generally, it follows that
any exact sequence $\ ...\ \to V_i\ \to V_{i+1}\ \to\ ...$ of $B_n$ modules defines an exact sequence
$\ ...\ \to V_{i+1}^t\ \to V_i^t\ \to\ ...$ of $B_n$ modules.
\end{remark}

\subsection{Conventions}\label{conventions:sec} While we are primarily interested in representations over the complex numbers $\C$, it will be useful to deal with our algebras in a more general setting. In the following, we assume $F_0$ to be a field of characteristic 0 which contains a primitive third and fourth root of unity. We define $R=F_0[\la_i^{\pm 1}, (\la_i-\la_j)^{-1}]$, for $1\leq i\neq j\leq 3$, where the $\la_i$ are viewed as variables. 
$F$ is assumed to be a field with $\la_i\in F$, $\la_i\neq\la_j$ for $i\neq j$. We will usually encounter $F$ as the quotient field $Q(R/\p)$ for a prime ideal $\p\subset R$. We may write $F\K_n$ if we want to stress the field over which we define the algebra; otherwise we just write $\K_n$.

Statements concerning regular modules usually mean that it holds for the given module $V$ or its dual module $V^t$. The discerning reader should replace statements dealing with non-semisimple regular modules as: For the given regular module $V$ or its dual module $V^t$ the following holds... .

\subsection{The algebra $\K_3$}\label{K3:sec}
We determine the structure of the algebra $\K_3$ in this section. These results are known to experts, see e.g. the review in \cite{CJ}, Chapter 2. They could also be derived e.g. from results in \cite{TW}.

\begin{theorem}\label{B3class} \cite{BM}, \cite{TW} The algebra $\K_3$ is semisimple over a field $F$ of characteristic 0 except if one of the following polynomials in the eigenvalues is equal to 0:
$$\la_i-\la_j,\quad \la_i^2-\la_i\la_j+\la_j^2\quad {\rm or}\quad \la_i^2+\la_j\la_k,$$ where $\{ i,j,k\} = \{ 1, 2, 3\}$.
In the semisimple case, we have three 1-dimensional representations $\{ \la_i\}$, three two-dimensional representations
$\{ \la_i\la_j\}$ and one three-dimensional representation $\{ \la_1\la_2\la_3\}$. Here the expression between the
brackets $\{ \ \}$ is the determinant of $\sigma_1$ in the given representation.
\end{theorem}

We also expect the following result to be known.

\begin{corollary}\label{K3quot} Replacing a regular representation $V$ by its conjugate $V^t$ if necessary, see Remark \ref{pathunique},
we obtain the following exact sequences, which do not split:

\noindent (a) If $\la_i=-\theta\la_j$, we have
\begin{equation}\label{B3}
0 \to \{\la_i\}\to\{\la_i\la_j\}\to \{\la_j\}\to 0.
\end{equation}
(b)  If $\la_1^2+\la_2\la_3=0$, we have
\begin{equation}\label{B32}
0\to \{\la_2\la_3\}\to \{ \la_1\la_2\la_3\}\to\{\la_1\}\to 0.
\end{equation}
\end{corollary}

$Proof.$ The first exact sequence is known from the study of Iwahori-Hecke algebras, see Section \ref{sec:regularmod}, (b) for more details. Statement (b) follows e.g. from Corollary \ref{6dimquot},(a). It can also be shown directly using the formulas
for braid matrices in Theorem \ref{AB2theorem} and Corollary \ref{nonsemisimp}, and the results in Section \ref{examplesnthree}, Example (1).

\begin{lemma}\label{coincidenceK3} Assume that the eigenvalues are mutually distinct, and two of the equations  $\la_i^2-\la_i\la_j+\la_j^2=0$ or $\la_i^2+\la_j\la_k=0$ in Theorem \ref{B3class} are satisfied. 
Then we have $\la_j=-\theta\la_k$, $\la_k=-\theta\la_i$  for a primitive 3rd root of unity $\theta$ and indices $\{ i,j,k\}=\{ 1,2,3\}$. In particular, all simple representations of $\K_3$ are 1-dimensional in this
case, except for $\{\la_i\la_k\}$.
\end{lemma}

$Proof.$ If $\la_i^2=-\la_j\la_k$ and $\la_j^2=-\la_i\la_k$, we can solve for $\la_k$ in both equations. It follows that $\la_i^2/\la_j=\la_j^2/\la_i$, which implies that $\la_i=\theta\la_j$ for some primitive third root of unity $\theta$ (as we assumed $\la_i\neq\la_j$).
Substituting this into  $\la_j^2=-\la_i\la_k$, we obtain  $\la_j=-\theta\la_k$ and $\la_k=-\theta\la_i$.  Similar proofs work for
any choice of the equations mentioned in the statement. We can rule out cases $\la_i=-\theta\la_j$ and $\la_j=-\theta^2\la_k$,
as this would imply $\la_i=\la_k$. This is excluded in our setting. 

Assume $\la_1=-\theta\la_2$ and $\la_2=-\theta\la_3$.
Then, by Corollary \ref{K3quot}, the representations $\{\la_1\la_2\}$ and $\{\la_2\la_3\}$ are not simple. Hence any subquotient module in these modules must be 1-dimensional.
Moreover, these equations also imply $\la_1^2+\la_2\la_3=0=\la_3^2+\la_1\la_2$. Hence $\{\la_1\la_2\la_3\}$ has a 1-dimensional submodule, and a 2-dimensional quotient $\{\la_2\la_3\}$ which is not simple. 
This shows that $\{\la_1\la_2\la_3\}$ has a composition series with one-dimensional factors.

\begin{remark} Using the explicit results in \cite{TW}, it is also possible to determine all simple modules of $\K_3$ when we have repeated eigenvalues. E.g. if $\la_1=\la_2=\la_3$,
it follows from \cite{TW}, Theorem 2.9 that we have exactly one simple module of dimensions 1, 2 and 3. We will not pursue these cases further in this paper.
\end{remark}

\subsection{Regular modules of $\K_4$}\label{genericirrep}  We want to carry out a similar analysis for the algebra $\K_4$. As a first step, we recall what structure it has
when it is semisimple.

\begin{theorem}\label{K4generic} \cite{MM} Let $F$ be a field of characteristic zero which contains a primitive third root of unity $\theta$.
The algebra $\K_4$, defined over $F$,
 has the same dimension and, for generic eigenvalues,
the same decomposition into a direct sum of full matrix rings  as the group algebra of the complex reflection group $G_{25}$.   A full list of regular modules, up to permutation of the eigenvalues, is listed in Table 1.
\end{theorem}
The regular modules of the algebra $\K_4$ are already determined by
the determinant det$(\sigma_1)$ of a generator of the braid group, except for the 9-dimensional representations, which
also depend on the choice of a primitive third root of unity $\theta$. We will refer to a representation
 by $\{ {\rm det}(\sigma_1)\}$, where we add the suffix $\theta$ to the  nine-dimensional representation $\{ \la_1^3\la_2^3\la_3^3\}_\theta$.
  We list the scalar by which the central element
 $\Delta_4^2$ acts in the third  column below. The last column tells us the decomposition of the simple $\K_4$
module into a direct sum of $\K_3$ modules. Weights will be defined in Section \ref{sec:weightspaths}. $Convention$: $V=X+Y$ means that $V$ has a composition series with factors $X$ and $Y$.

\[
\begin{tabular}{|c|c|c|c|c|} dim & det($\sigma_1$)&$\Delta_4^2$&weights& $\K_3$ rep\\
\hline
& & & &\\
9 & $\la_1^3\la_2^3\la_3^3$&$ \theta\la_1^4\la_2^4\la_3^4
$& $(i,j),\ 1\leq i,j\leq 3$& see $R9$\\
8&$\la_1^4\la_2^2\la_3^2$&
$\la_1^6\la_2^3\la_3^3$&  (1,1): 2, $(i,j), i\neq j$: 1 & see $R8$\\
6&$\la_1^3\la_2^2\la_3$&
$\la_1^6\la_2^4\la_3^2$& $(1,i),(i,1)$ and $(2,2)$& see $R6$ \\
3 & $\la_1\la_2\la_3$ & $\la_1^4\la_2^4\la_3^4$& $(i,i)$&$\{\la_1\la_2\la_3\}$\\
3 & $\la_1^2\la_2$ & $\la_1^8\la_2^4$ & $(1,1),(1,2),(2,1)$&$\{\la_1\la_2\} + \{\la_1\}$\\
2 &  $\la_1\la_2$ & $\la_1^6\la_2^6$ & $(1,1),(2,2)$&$\{\la_1\la_2\}$\\
1 & $\la_1$ & $\la_1^{12}$ & $(1,1)$&$\{\la_1\}$\\
\end{tabular}
\]
\vskip .2cm
$R9:$\quad $\{\la_1\la_2\la_3\} +\{\la_1\la_2\}  +\{\la_1\la_3\}+\{\la_2\la_3\}$,

$R8:$\quad  $\{\la_1\la_2\la_3\}  + \{\la_1\la_2\}  +\{\la_1\la_3\} + \{\la_1\}$,

$R6:$\quad  $\{\la_1\la_2\la_3\}  +\{\la_1\la_2\}  +\{\la_1\}$.

$${\rm Table\ 1}$$
\begin{remark}\label{prelimrem}  (a) Let $\Ga = \{ \ga\}$ be a labeling set for the isomorphism classes of simple modules of $\K_4$, and let $\la_\ga$ be the scalar
via which $\Delta_4^2$ acts on the regular module $V_\ga$. If $\la_\ga\neq\la_{\ga'}$ for all $\ga'\neq\ga$, the eigenprojection $z_\ga$ of $\Delta_4^2$ for this scalar is a central idempotent of $\K_4$. It corresponds to the simple representation $V_\ga$ in the semisimple case. Unfortunately, we are not aware of a general result which would imply that the generic representation $V_\ga$ will be simple whenever $z_\ga$ is well-defined, even though we will be able to prove this later, with some effort. 

(b) Assuming the result in (a), we could narrow down the values for the eigenvalues for which $\K_4$ is not semisimple.
It suffices to just set equal the eigenvalues of $\Delta_4^2$ for two generic representations.
One could thus show that $\K_4$ is not semisimple only if one of the following expressions is equal to zero: $\la_i^{12}-\la_j^{12}$,
 $\la_i^6\pm \la_j^3\la_k^3$ or $\la_i^3-\la_j^2\la_k$, where $\{ i,j,k\}=\{ 1,2,3\}$.
Compare this with the precise result, see Theorem \ref{semisimple:thm}.
\end{remark}

\subsection{Weights and paths}\label{sec:weightspaths} Unless noted otherwise, we assume that we only deal with representations of $\K_4$ in which the images of the braid generators are diagonalizable. While regular representations are well-defined for any choice of eigenvalues, 
it is hard to see for which values of the eigenvalues they are not simple.
As we are primarily interested in the case where the images of the braid generators are diagonalizable, we can construct different bases
which are more useful in that aspect. Observe that
the images of $\sigma_1$ and $\sigma_3$ form a commutative subalgebra of the image of $B_4$ in any representation
of $B_4$. This motivates the following definition.

\begin{definition}\label{weights} Let $W$ be a representation of $B_4$. We call any common eigenvector of $\sigma_1$
and $\sigma_3$ a {\it weight vector}. If $\sigma_1$ has eigenvalues $\la_1,\la_2,\ldots,\la_d$, we say that a vector $v\in W$ has
weight $(i,j)$ if $\sigma_1$ and $\sigma_3$ act on it via multiplication by $\la_i$ and $\la_j$ respectively.  The dimension of the corresponding weight space $V(i,j)$
is called the multiplicity of the weight. 
\end{definition}
\begin{remark}\label{weightsymmetry}
It follows from the definitions and \ref{braidpermutation} that $\Delta_4$ maps the weight space for the weight $(i,j)$ onto the weight space for the weight $(j,i)$.
\end{remark}
The weights and their multiplicities are listed in the second to last column in  Table 1. We see that the weights have multiplicity at most 1
in all regular modules, except for $V=\{\la_r^4\la_s^2\la_t^2\}$, where the weight $(r,r)$ has multiplicity 2. This will allow us to construct convenient bases for regular modules even in the case of multiplicity $>1$, see Section \ref{sec:weights}.

We will also use another well-known type of bases, {\it path bases}. They can be defined for any sequence of finite-dimensional semisimple algebras $\A_n$ for which
the restriction of a simple $\A_n$ module to $\A_{n-1}$ decomposes into a direct sum of mutually nonisomorphic $\A_{n-1}$ modules. A standard example
would be $\A_n=\C S_n$, the group algebra of the symmetric group $S_n$. We will only give the definitions for the special case in this paper.

Let us assume that the algebras $\K_2\subset \K_3\subset \K_4$ are semisimple. It follows from
the restriction rules stated in Table 1  that any simple 
$\K_n$ module $V_n$ decomposes into a direct sum of mutually nonisomorphic $\K_{n-1}$ modules for $n=3,4$ in the generic case. 
Hence, if $\ga_i$ labels
a simple $\K_i$ module, and $z_{\ga_i}$ is the corresponding central idempotent in $\K_i$,  it follows  that 
dim $z_{\ga_2}z_{\ga_3}z_{\ga_4}V\leq 1$ for any regular $\K_4$ module $V$.
\begin{definition}\label{def:paths}
If the product $z_{\ga_2}z_{\ga_3}z_{\ga_4}$ is nonzero, we call $t=(\ga_i)$
a $path$, the idempotent $p_t=z_{\ga_2}z_{\ga_3}z_{\ga_4}$ a {\it path idempotent} and any nonzero $v_t\in p_tV$ a {\it path  vector} for $V$.  We also
define $t(i)=\ga_i$.
\end{definition}

\subsection{Bases of path vectors or weight vectors} We shall use bases of path or weight vectors to determine when a regular module is simple.
We list a number of simple properties for such bases in the following general setting. Let $\A$ be an algebra over a field $F$, and let $V$ be a finite-dimensional $\A$-module.
We assume we have a family $(P_t)_{t\in I}\subset \End(V)$ of rank 1 idempotents in the image of $\A$ such that $P_tP_s=0$ for $t\neq s$ and a basis $(v_t)$
with $v_t\in P_tV$. We define a relation
\begin{equation}\label{def:relation}
s\leftarrow t\ \Leftrightarrow\ P_s\A v_t\neq 0,
\end{equation}
i.e. there exists an element $a\in \A$ such that its $st$ matrix coefficient with respect to our basis is not equal to zero. We also denote
by $\leftarrow$ the transitive closure of this relation, and we write $s\leftrightarrow t$ if $s\leftarrow t$ and $t\leftarrow s$. We then have
the following lemma, whose proof follows from the definitions and basic linear algebra.

\begin{lemma}\label{matrixcoeff}  Let $0=V_0\subset V_1\subset\ ...\ \subset V_m=V$ be a composition series with simple quotients $V_i/V_{i-1}$.
Then we obtain bases for the factors $V_i/V_{i-1}$ from the vectors $v_t$ labeled by indices $t$ in an equivalence class of the relation $\leftrightarrow$.
In particular, $V$ is simple if all indices are equivalent to each other.
\end{lemma}

\medskip
We collect a few basic facts about path bases. More details can e.g. be found in \cite{MW2} for more general braid representations. An explicit
example for the matrix blocks mentioned in the following lemma appears in Lemma \ref{8dimcomplem}.

\begin{lemma}\label{pbaseprop} Let $V$ be a $\K_4$ module, and let $S_i$ be the matrix representing $\sigma_i$ for $i=2,3$ with respect to  a path basis $(v_t)$.

(a) Let $t=(\ga_j)$ be a path. Then $S_iv_t$ is a linear combination of paths $v_s$ with $s\in [t;i]$, the set of all paths $s=(\nu_j)$ such that $\nu_j=\ga_j$ for all $j\neq i$.
So, $S_i$ is a direct sum of matrix blocks acting on the subspaces spanned by vectors $v_s$ with  $s\in [t;i]$.

(b) Assume one of the eigenprojections $E$ of $S_i$ has rank 1 in the block labeled by the paths in $[t;i]$. 
Then  $s\leftrightarrow s'$ for $s,s'\in [t;i]$ if $d_sd_{s'}\neq 0$ for the diagonal entries $d_s$, $d_{s'}$ of $E$.
\end{lemma} 

$Proof.$ Statement (a) follows from the fact that $S_i$ commutes with $z_{t(j)}$ for all $j\neq i$. For part (b), let $P_s$ be the projection on the span of $v_s$.
If $s,s'\in [t;i]$, we have $P_sEP_{s'}EP_s=d_sd_{s'}P_s\neq 0$, using rank $E=1$. This implies $s\leftarrow t$ and $t\leftarrow s$.

\subsection{Regular modules}\label{sec:regularmod} Here we list explicit regular modules which we will use to study the representation theory of the algebras $\K_3$ and  $\K_4$.

(a) For modules of $\K_n$, $n=3,4$, in which $\sigma_i$ acts with only two distinct eigenvalues, we obtain a representation of the Iwahori-Hecke algebra $\H_n$ of type $A_{n-1}$.  A class of special  modules, called Specht modules, were defined in \cite{DJ1} for all Hecke algebras $\H_n$ which satisfy the conditions of regular modules.  Composition series for Specht modules have been determined in general, see e.g. \cite{Mathas}, \cite{GJ}. We will use the Specht modules of $\H_3$ and $\H_4$ as regular modules for the representations labeled by  $\{\la_i\}$, $\{\la_i\la_j\}$ (for both $\K_3$ and $\K_4$) and $\{\la_i^2\la_j\}$.  We can deduce  the results already mentioned in Section \ref{K3:sec} from \cite{DJ1} and \cite{DJ2} or the aforementioned lecture notes, as well as the following
exact sequence when $\la_i^2+\la_j^2=0$:
\begin{equation}\label{Hecke4}
0\to \{\la_i\}\ \to\ \{\la_i^2\la_j\}\ \to\ \{\la_j^2\la_i\}\to \{\la_j\}\ \to\ 0.
\end{equation}
In particular, we obtain a simple two-dimensional representation $\{\la_i\la_j\}^*$. Unlike the $\K_4$ module $\{\la_i\la_j\}$, it has the weights $(i,j)$ and $(j,i)$.

(b) The  $\K_4$ modules $\{\la_i\}$, $\{\la_i\la_j\}$ and $\{\la_1\la_2\la_3\}$ are  obtained from the modules of $\K_3$ denoted by the same symbol by mapping $\sigma_3$ to the matrix representing $\sigma_1$.
 In particular, the results   in Section \ref{K3:sec} also apply for these $\K_4$ modules.

(c)  By \cite{MM}, Proposition 7.3, each irreducible representation of the Hecke algebra $\K_4$  does appear in some pre-$W$-graph. Explicit matrix representations  for the 6, 8 and 9-dimensional 
representations were determined in \cite{BM}, \cite{MM}, based on this; see \cite{MWa}, Section 4.4, Table 4 for explicit matrices.
 In all these representations, the entries of the matrices representing
the braid generators are in the ring $\Z[\theta][\la_1,\la_2, \la_3]$, where $\theta$ is a primitive third root of unity, i.e. they are regular modules in our definition.
An illustration for the paths for semsimple $\{\la_i^3\la_j^2\la_k\}$ and $\{\la_i^4\la_j^2\la_k^2\}$ are given by the  paths in the Bratteli diagram in \cite{MW}, Section 2.1 from 0 to $3\La_1$ and $\La_1+\La_2$ respectively. 
See the proof of Lemma \ref{9dimcomplem} for paths for the 9-dimensional representation.

\subsection{Permutation automorphisms} We consider the algebras $\K_n$, $n\leq 5$, defined over the polynomials
$R=F_0[\la_1,\la_2,\la_3]$. Obviously, any permutation $\pi\in S_3$ induces an automorphism of $R$ via permuting the eigenvalues.
We can extend this to an automorphism of $\K_n$, also denoted by $\pi$ via the map
$$\pi: \sum_{\beta\in B_n} a_\beta\beta\ \mapsto \  \sum_{\beta\in B_n} \pi(a_\beta) \beta;$$
it is clear that this map factors over the defining ideal for $\K_n$. We list some simple properties
of these permutation automorphisms.

\begin{lemma}\label{permautlem} Let $\pi\in S_3$ also denote the corresponding automorphism of $R=F[\la_1,\la_2,\la_3]$ and of $\K_n$. Then we have

(a) $\pi$ induces a permutation of the simple components of $\K_n$, $n\leq 5$. If $\{\la_1^{e_1}\la_2^{e_2}\la_3^{e_3}\}$ denotes a simple representation of $\K_n$, $n\leq 4$, then $\pi$ maps it to the simple representation denoted by
 $\{\pi(\la_1^{e_1}\la_2^{e_2}\la_3^{e_3})\}$.

(b) $\pi$ induces a permutation of path idempotents.

(c)  If $p_t$ is a path idempotent, and $P_{ki}$ the eigenprojection of $S_k$ for the eigenvalue $\la_i$, then
$$\pi(p_t)P_{k\pi(i)}\pi(p_t)\ =\pi(\al)\pi(p_t),$$
where $\al\in R$ is defined by $p_tP_{ki}p_t=\al p_t$.
\end{lemma}

\subsection{$BMW$ algebras}\label{BMW:sec} It has been independently shown in \cite{BBMW} and \cite{Mu} that one obtains a finite dimensional quotient algebra $\BM_n$ of $FB_n$ for $all$ $n$
if we add to the cubic relation for $\K_n$ one more relation, see the quoted papers for details. If $P_1(\la_3)$ denotes the eigenprojection of $\sigma_1$ for the eigenvalue $\la_3$,
the relation could be expressed as 
$${\la_3^{\pm 1}(\la_1-\la_3)(\la_2-\la_3)}\ P_1(\la_3)\sigma_2^{\pm 1}P_1(\la_3)\ =\ {\la_3(\la_1+\la_2)}P_1(\la_3),$$
with matching signs in the exponents.
It was shown in \cite{BBMW} and \cite{Mu} that the simple components of $\BM_n$ are labeled by Young diagrams $\mu$
with $|\mu|=n-2k$ boxes, $0\leq k\leq \lfloor n/2\rfloor$, in the generic semisimple case. 
If $V_{n,\mu}$ denotes a simple $\BM_n$ module, it decomposes as a $\BM_{n-1}$
module into the direct sum of modules $V_{n-1,\ga}$, where $\ga$ runs through the diagrams obtained by removing or, if $|\mu|\leq n-2$,
adding a box from/to $\mu$. One deduces from this that $\BM_3\cong \{\la_1\la_2\la_3\}\oplus \{ \la_1\la_2\}\oplus \{\la_1\}\oplus \{ \la_2\}$,
and that 
\begin{equation}\label{BMW4}
\BM_4\cong \{\la_1^3\la_2^2\la_3\}\oplus \{\la_1^2\la_2^3\la_3\}\oplus \{\la_1\la_2\la_3\}\oplus \H_4(\la_1,\la_2),
\end{equation}
where $\H_4(\la_1,\la_2)$ is the direct sum of all simple representations of $\K_4$  in which the eigenvalue $\la_3$ does not appear.
Here the eigenvalue $\la_3$ corresponds to the representation of $\BM_2$ labeled by the empty diagram $\emptyset$,
$\la_1$ corresponds to $[2]$ and $\la_2$ corresponds to the diagram $[1,1]$. One checks that
for $\BM_3$, $[1]$  labels the representation $\{\la_1\la_2\la_3\}$, and for $\BM_4$, 
$[2]$ labels $\{\la_1^3\la_2^2\la_3\}$, $[1,1]$  labels $\{\la_1^2\la_2^3\la_3\}$ and  $\emptyset$  labels $\{\la_1\la_2\la_3\}$.

The eigenvalues of the braid generators of the algebras $\BM_n=\BM_n(r,q)$ in the original papers were normalized such that
 $\la_1=q$, $\la_2=-q^{-1}$ and $\la_3=r^{-1}$. If $r=\pm q^k$ for an integer $k$, the algebras $\BM_n$ may no longer 
be semisimple. We have the following result, see \cite{WBMW}, Corollary 5.6.

\begin{theorem}\label{BMWquot} If $r=\pm q^k$ for some integer $k$, and $q$ is not a root of unity, the algebra $\BM_n$ has a semisimple quotient,
whose simple components are labeled by Young diagrams $\mu$ with $n-2k$ boxes, $0\leq k\leq \lfloor n/2\rfloor$, where now $\mu$ can only 
be taken from a restricted subset $\Lambda(\pm , k)$ of Young diagrams. The dimensions of its simple representations can be calculated
by the same restriction rule as before, but with only Young diagrams from $\Lambda(\pm,k)$ allowed.
The set  $\Lambda(+,N-1)$ consists of all Young diagrams with $\leq N$ boxes in the first two columns, while the set 
 $\Lambda(-,N+1)$ consists of all Young diagrams with $\leq  N$ boxes in the first column.
\end{theorem}
\begin{remark} The set $\Lambda(+,N-1)$ labels the irreducible representations of $O(N)$ which appear in tensor powers of
its vector representation $\C^N$, while  the set $\Lambda(-,N+1)$ labels the irreducible representations of $Sp(N)$
for $N$ even. In particular, we can identify the eigenprojections for the eigenvalues $\la_1,\la_2$ and $\la_3$ in $\K_2$
with the irreducible representations of $O(N)$ and $Sp(N)$ labeled by the Young diagrams $[2]$, $[1,1]$ and $\emptyset$.
It follows from the definitions that a representation of $\BM_4$ labeled by the Young diagram $\mu$ has weight $(i,j)$
with multiplicity $m$
just corresponds to the fact that the $O(N)$ or $Sp(N)$ representation labeled by $\mu$ appears with multiplicity $m$ in
the tensor product of the representations corresponding to $\la_i$ and $\la_j$.
\end{remark}

\begin{corollary}\label{6dimquot} (a) If $\la_2^2=-\la_1\la_3$, we obtain a two-dimensional quotient of the $\BM_3$ representation $\{\la_1\la_2\la_3\}$
which is isomorphic to $\{\la_1\la_3\}$.
\vskip .2cm
(b) We obtain nontrivial quotients of the 6-dimensional representation $\{\la_1^3\la_2^2\la_3\}$ as follows:
\vskip .2cm
 (1) If $\la_1=-\la_3$,  we obtain a 3-dimensional quotient  $\{\la_2|\la_1\la_3\}^*$  in $\{\la_1^3\la_2^2\la_3\}$. 
It has the weights $(22),(13)$ and $(31)$, and it is isomorphic to $\{\la_1\la_2\la_3\}$ as a $\K_3$ module. The submodule is isomorphic to the generic representation $\{\la_1^2\la_2\}$.
\vskip .2cm
(2) If $\la_2=-\la_3$, we obtain a 4-dimensional quotient $\{\la_1^2\la_2\la_3\}$  with  weights
$(12),(21), (13)$ and $(31)$. It is isomorphic to $\{\la_1\la_2\la_3\}+\{\la_1\}$ as a $\K_3$ module. The submodule is the generic $K_4$ module $\{\la_1\la_2\}$.
\vskip .2cm
 (3) If $\la_2^2=-\la_1\la_3$, we have a quotient and a submodule which are isomorphic to the generic Hecke algebra representations 
$\{\la_1^2\la_3\}$ and $\{\la_1\la_2^2\}$.
\vskip .2cm
(4) If $\la_1^3=\la_2^2\la_3$, we obtain a 5-dimensional quotient $\{\la_1^2\la_2^2\la_3\}$ with the  weights
$(22),(12),$ $(21),(13)$ and $(31)$.  It is isomorphic to $\{\la_1\la_2\la_3\}+\{\la_1\la_2\}$ as a $\K_3$ module. 
\vskip .2cm
(c) All subquotient modules described in (a) and (b) are simple, as long as the eigenvalues do not satisfy an additional relation.
\end{corollary}

$Proof.$ The statements do not change if we multiply the images of the braid generators by a constant. It is therefore enough to prove the statements with the eigenvalues $\la_1=q$, $\la_2=-q^{-1}$ and $\la_3=r^{-1}$. 

We start with proving the statements in part (b). Observe that $\la_2=-\la_3$ implies $r^{-1}=q^{-1}$. Hence 
we are in the case $(+,1)$ corresponding to $O(2)$, where only diagrams are allowed with $\leq 2$ boxes in the first two columns. 
It follows that the $\BM_4$ representations labeled by $[2]$ and $[1,1]$ have dimensions 4 and 3 respectively. One derives from the tensor product rules of $O(2)$ that the representation $V_{[2]}$ appears in the tensor products $V_\emptyset \otimes V_{[2]}$
and $V_{[1,1]}\otimes V_{[2]}$. The same applies if we change the order of the tensor factors. Hence the image of $\BM_4$
in the representation labeled by $[2]$ has the weights $(1,2), (2,1), (1,3)$ and $(3,1)$. This proves (2), using the fact that
the kernel of this representation must have weights $(1,1)$ and $(2,2)$. Similarly, we deduce from the tensor product rules
of $O(2)$ that the representation $V_{[1,1]}$ appears in $V_{[2]}\otimes V_{[2]}$, and its tensor product with the trivial
representation $V_\emptyset$. This shows the 3-dimensional quotient of $\{\la_2^3\la_1^2\la_3\}$ has the weights $(1,1)$, $(2,3)$ and $(3,2)$. Applying the permutation automorphism for the transposition (12), we obtain statement (1). 

Observe that $r=-q^3$ implies that $\la_2^2=-\la_1\la_3$ for our choice of eigenvalues. So this corresponds to the case of $Sp(2)$,
where only diagrams with one row are allowed. It follows that only the representations corresponding to $\{\la_1\}$ and $\{\la_3\}$
are allowed in the image of $\BM_2$ and the representation $\{\la_1\la_2\la_3\}$ has a quotient isomorphic to $\{\la_1\la_3\}$.
This proves statement (a). We leave it to the reader to check that we get a 3-dimensional quotient
for the $\BM_4$ representation labeled by $[2]$ with weights $(1,1)$, $(1,3)$ and $(3,1)$ which corresponds to
the usual representation $\{\la_1^2\la_3\}$. The kernel has the same dimension and weights as the usual $\{\la_2^2\la_1\}$.

One similarly links the case $\la_2^3=\la_1^2\la_3$ to the $\BM_4$ representation with $r=-q^5$, which corresponds to 
$Sp(4)$. As the diagram $[1,1,1]$ is no longer allowed, the $\BM_4$ representation labeled by $[1,1]$ only has dimension 5.
It follows that the $\K_4$ representation $\{\la_2^3\la_1^2\la_3\}$ has a 5-dimensional quotient with weights $(1,1)$, $(2,1), (1,2)$, $(2,3)$ and $(3,2)$. Applying the permutation automorphism for the transposition $(12)$ gives us (4).

It was shown in the proof of \cite{WBMW}, Corollary 5.6 that the quotient modules are simple unless $q$ is a root of unity. One can check that this 
is equivalent with the statement about additional relations. This statement can also be directly proved with the methods in this paper.

\section{Weight bases}\label{weight base:sec}
We  assume for the whole section that we only deal with representations of $\K_4$ in which the images of the braid generators are diagonalizable.

\subsection{Weight bases and semisimplicity}\label{sec:weights} Recall the definition of weights in Section \ref{sec:weightspaths}.
We see in Table 1 that the weights have multiplicity at most 1
in all regular modules, except for $V=\{\la_r^4\la_s^2\la_t^2\}$, where the weight $(r,r)$ has multiplicity 2. 
We also want to define canonical weight vectors in the latter case.
In the following, let $S_i$ be the matrix representing $\sigma_i$ in a fixed regular representation, and let
\begin{equation}\label{Piidef}
P_{kr}=\prod_{s\neq r} (S_k-\la_s1)\quad {\rm and }\quad B(i,j)=P_{1i}P_{3j}.
\end{equation}
It should be clear that
$P_{kr}=(\la_r-\la_s)(\la_r-\la_t)P_k(\la_r),$
where $P_k(\la_r)$  is the eigenprojection  of $S_k$ for its eigenvalue $\la_r$ in $V$, and that
 the image of 
$$B(i,j)=P_{1i}P_{3j}$$
is indeed the weight space of $(i,j)$.  We will first need the following result, where the second statement of (b) can be skipped for the first reading.

\begin{lemma}\label{8dimcalc} (a) The element $\Delta_4$ maps the weight space $V(i,j)$ onto the weight space $V(j,i)$ for any regular $\K_4$ module $V$.

(b) If $V=\{\la_1^4\la_2^2\la_3^2\}$, $\Delta_4$ acts on the 2-dimensional weight space $V(1,1)$ with eigenvalues $\pm\la_1^3(\la_2^3\la_3^3)^{1/2}$.
\end{lemma}

$Proof.$ Part (a) is a direct consequence of Eq \ref{braidpermutation}. For the proof of (b), we use the explicit 8-dimensional representation in \cite{MWa}, with $a=\la_r$, $b=\la_s$, $c=\la_t$.
One calculates that  the matrix $B$ which represents $\Delta_4P_1(\la_1)P_3(\la_1)$ has
nonzero entries only in the third and fifth row. The  $2\times 2$ submatrix for these rows and columns is given by
$$a^3bc\ \left[\begin{matrix} c&c(b-c)\cr 1&-c\end{matrix}\right].$$ 
This is all one needs to check that the characteristic polynomial of $B$ is equal to $\la^6(\la^2-a^6b^3c^3)$.

\begin{definition}\label{pathprojections} Let $V$ be a regular $\K_4$ module. If $(r,s)$ is a weight of $V$ with multiplicity 1, we call $P_1(\la_r)P_3(\la_r)$ its canonical weight projection.
If the weight $(r,r)$ has multiplicity 2, its canonical weight projections are defined to be the eigenprojections of $\Delta_4P_1(\la_r)P_3(\la_r)$ for its non-zero eigenvalues, after possibly having to adjoin a square root of $a^6b^3c^3$ to the ground ring.
A {\it weight basis} for $V$ is obtained by picking a nonzero vector in the image of each weight projection.
\end{definition}


\subsection{Conditions for simplicity of 6-dimensional representations $\{ \la_i^3\la_j^2\la_k\}$} We will use weight projections
to construct a full system of matrix units for regular representations whenever the eigenvalues satisfy certain equations.
We use the notations from \ref{Piidef}.
\begin{lemma}\label{6dimweightbasis}
The 6-dimensional module $\{\la_1^3\la_2^2\la_3\}$  is simple except possibly if $\la_1=-\la_3$, $\la_2=-\la_3$,  $\la_1^3=\la_2^2\la_3$
or $\la_2^2=-\la_1\la_3$.
\end{lemma}

$Proof.$ It follows from Definition \ref{pathprojections} 
that the module $\{\la_1^3\la_2^2\la_3\}$ has a basis of weight vectors whose weights are listed in
Table 1.
Using the explicit matrix representations in \cite{MWa} and,
preferably, some symbolic computing system, we calculate the quantity $d(i,j)$ from the identity
$$B(i,j)P_{23}B(i,j) =d(i,j)B(i,j).$$
The quantities $d(i,j)$ are the diagonal entries of $P_{23}$, up to a constant multiple involving factors of the form $\la_i-\la_j$.
It follows from \ref{braidpermutation} that conjugation by $\Delta_4$ fixes $P_{23}$ and permutes $B(i,j)$ 
with $B(j,i)$; this implies $d(i,j)=d(j,i)$.  We obtain the values
\begin{align}
d(1,1)\ &=\ (\la_1-\la_3)(\la_1+\la_3)(\la_2+\la_3)(\la_1^3-\la_2^2\la_3),\cr
d(2,2)\ &=\ \la_1(\la_2-\la_3)^2(\la_2+\la_3)(\la_2^2+\la_1\la_3),\cr
d(1,2)\ =\ d(2,1)\ &=\ -\la_2(\la_1-\la_3)^2(\la_1+\la_3)(\la_2^2+\la_1\la_3),\cr
d(1,3)\ =\ d(3,1)\ &=\ (\la_1-\la_2)(\la_1-\la_3)(\la_1+\la_2)^2\la_3^2.\cr
\end{align}
It follows from Lemma \ref{matrixcoeff} that the module $V=\{ \la_1^3\la_2^2\la_3\}$ is simple if
all diagonal entries of the rank 1 matrix $P_{23}$ are nonzero.
This proves the claim, except that it could not be semisimple also if $\la_1=-\la_2$.
In this case, we still have a simple submodule $U$ of dimension at least 4, containing the
weights (11), (12), (21) and (2,2), as the corresponding diagonal entries of $P_{23}$ are nonzero.
As $P_{12}$ equals the sum of the two weight projections for (2,1) and (2,2) in our representation,
it follows that $U$ must contain all $\K_3$ submodules in which $P_{12}$ acts nonzero, i.e. $\{\la_1\la_2\la_3\}$
and $\{\la_1\la_2\}$. These are simple for $\la_1=-\la_2$, so $U$ must have dimension at least 5. It follows that it must have one of the weights of (13) or (31), and therefore,
 by Lemma \ref{8dimcalc}(a), both weights. This shows that $U=V$.

\subsection{Conditions for simplicity of 8-dimensional and 9-dimensional representations}

\begin{lemma}\label{9and8dim} Assume that $\{ r,s,t\}=\{ 1,2,3\}$ and that $\la_r\neq \la_s$ for $r\neq s$.

(a) The representation $\{\la_1^3\la_2^3\la_3^3\}_\theta$ is simple except if $\la_r=-\theta\la_s$ or $\la_r\la_s=\theta\la_t^2$.

(b) The representation $\{\la_1^4\la_2^2\la_3^2\}$ is simple, except possibly if $\la_2=-\la_3$, $\la_3^3=\la_1^2\la_2$, $\la_2^3=\la_1^2\la_3$ or $\la_1^6=\la_2^3\la_3^3$.


\end{lemma}

$Proof.$ The proof is similar to the one of Lemma \ref{6dimweightbasis}, using Lemma \ref{matrixcoeff}. Let us do this first for the module $\{\la_1^3\la_2^3\la_3^3\}_\theta$. We define the scalar 
$\al_{ki,lj}$ by
$$B(k,i)S_2B(l,j)S_2B(k,i)=\al_{ki,lj}B(k,i).$$
It follows from an explicit calculation using the representation in \cite{MWa} that
$$\al_{11,12}\ =\  -\theta\la_1\la_2^2\la_3^2(\la_3+\theta\la_2)(\la_2+\theta\la_1)(\la_1\la_2-\theta\la_3^2)(\la_1\la_3-\theta\la_2^2)/(\la_1-\la_2).$$
Applying the permutation automorphism for the transposition $(23)$, see  Lemma \ref{permautlem}, we obtain $\al_{11,13}=(23)(\al_{11,12})$.
It follows by inspection that $\al_{11,12}$ and $\al_{11,13}$ are not equal to zero if the eigenvalues satisfy the conditions in the statement. Hence the weights $(1,j)$, $1\leq j\leq 3$ belong
to the same equivalence class of indices, see \ref{def:relation}. Applying an automorphism from Lemma \ref{permautlem} for a permutation $\pi$ with $\pi(1)=i$, $i=2,3$, we obtain the same statement for the weights $(i,j)$, $1\leq j\leq 3$.  
As $(i,j)$ and $(j,i)$ are equivalent by Remark 
\ref{weightsymmetry} for $1\leq i,j\leq 3$, it follows that all weights $(i,j)$ are equivalent as long as $\pi(\al_{11,12})\neq 0$ for all $\pi\in S_3$.
The claim now follows from Lemma \ref{matrixcoeff}.

We can proceed similarly for the representations $\{\la_i^4\la_j^2\la_k^2\}$. By symmetry, it suffices to do this for $V=\{\la_1^4\la_2^2\la_3^2\}$ as follows. 
We define $\al_{ij,kl}$ as before. We first calculate by brute force
$$\al_{21,23}=\la_1^3\la_3(\la_3^3-\la_1^2\la_2)(\la_2+\la_3)(\la_1-\la_2)(\la_2-\la_3).$$
As the permutation $(23): \la_2\leftrightarrow\la_3$ gives us an equivalent representation, we can argue as in the last paragraph that
$\al_{31,32}=(23)(\al_{21,23})$ is nonzero whenever $\la_i^3\neq \la_1^2\la_j$, where $\{ i,j\}=\{ 2,3\}$, and $\la_2\neq -\la_3$. Hence $(2,1)\leftrightarrow (2,3)$ and  $(3,1)\leftrightarrow (3,2)$, with $\leftrightarrow$ as defined before Lemma \ref{matrixcoeff}.
We can show $(i,j)\leftrightarrow (j,i)$ as in the previous paragraph. This already implies that our representation  $\{\la_1^4\la_2^2\la_3^2\}$ has a simple subquotient module containing the span $W$ of all vectors with weights $(i,j)$, $i\neq j$.
Assume $W$ is contained in a proper submodule $U$ of $V$. Then $V/W$ and hence also $V/U$ would only have the weight $(1,1)$. This entails that $\Delta_4^2$ would act on $V/U$ via the scalar $\la_1^{12}=\la_1^6\la_2^3\la_3^3$, which implies  $\la_1^6=\la_2^3\la_3^3$.

\subsection{Semisimplicity} We can now  almost prove one direction of one of the main results of this paper.
\begin{proposition}\label{semisimplevalues:prop} Assume $\la_i\neq \la_j$ for $i\neq j$.
A regular module $V$ is simple unless the eigenvalues satisfy one of the polynomials listed next to it in Table 2 below, or, possibly, if  they satisfy $\la_j+\la_k=0$ or $\la_i^2=\la_j\la_k$ for the module $\{\la_i^4\la_j^2\la_k^2\}$.
In particular, $\K_4$ is semisimple if the eigenvalues do not satisfy any of the polynomials listed here.
\end{proposition}

$Proof.$ It suffices to prove that all regular modules are simple. Such modules with different labels also have different weights, see Table 1. Hence they are not isomorphic, which will imply semisimplicity of $\K_4$.  For the Iwahori-Hecke algebra modules in $\K_4$, i.e. modules in which $\sigma_1$ satisfies a quadratic equation,
it is well-known that they are not simple if and only if  $(-\la_i/\la_j)$ is  an $n$-th root of unity, $n=2,3,4$, see the discussion in Sections  \ref{K3:sec} and \ref{sec:regularmod}. The module $\{\la_1\la_2\la_3\}$
is not simple if and only if $\la_i^2+\la_j\la_k=0$, by Theorem \ref{B3class} and Corollary \ref{K3quot}. It follows from Lemmas \ref{6dimweightbasis} and \ref{9and8dim} that also the 6-, 8- and 9-dimensional modules
have to be simple if the eigenvalues do not satisfy one of the equations in the statement. This is summarized in Table 2. It will be shown in Theorem \ref{semisimple:thm} that the module $V$ is indeed not semisimple if the eigenvalues
are in the locus of one of the ideals listed in Table 2.
\[
\begin{tabular}{|c|c|c|c|} regular module $V$ & not semisimple\\
\hline
& \\
 $\{\la_1^3\la_2^3\la_3^3\}_\theta$&$(\la_i+\theta\la_j)$, $(\la_i^2-\theta\la_j\la_k)$\\
$\{\la_1^4\la_2^2\la_3^2\}$&
$(\la_2^3-\la_1^2\la_3)$, $(\la_3^3-\la_1^2\la_2)$, $(\la_1^2-\theta^{\pm 1}\la_2\la_3)$\\
$\{\la_1^3\la_2^2\la_3\}$&
$(\la_1+\la_3)$, $(\la_2+\la_3)$, $(\la_2^2+\la_1\la_3)$, $(\la_1^3-\la_2^2\la_3)$ \\
$\{ \la_1\la_2\la_3\}$ & $(\la_i^2+\la_j\la_k)$\\
$\{\la_1^2\la_2\}$ & $(\la_1\pm \sqrt{-1}\la_2)$\\
$\{\la_1\la_2\}$ & $(\la_1+\theta\la_2)$\\
\end{tabular}
\]
\vskip .2cm
\centerline {Table 2}

\subsection{Existence of path bases} Path bases do not exist in general. The following lemma will show that we can use them when $\K_4$ is not semisimple.

\begin{lemma}\label{pathexistence} Let $V$ be a regular module, and let $\p$ be one of the prime ideals listed next to it in Table 2.
Then $V$ has a path basis over the quotient field $Q(R/\p)$.
\end{lemma}

$Proof.$ As we assumed $\sigma_1=\Delta_2$ to be diagonalizable, it is sufficient to check that $\Delta_3^2$ acts via different eigenvalues on the $\K_3$ submodules of $V$.
This is fairly straightforward to check, using Table 1.
E.g if $V=\{\la_1^3\la_2^3\la_3^3\}_\theta$, such eigenvalues would only coincide if $\la_i^3=\la_j^3$ or $\la_i^2=-\la_j\la_k$, where $\{ i,j,k\}=\{ 1,2,3\}$.
This is not possible in $R/\p$ for $\p=(\la_i+\theta\la_j)$ or  $\p=(\la_i^2-\theta\la_j\la_k)$, as none of the former polynomials is in the later ideals. The same strategy works for 
$V=\{\la_1^4\la_2^2\la_3^2\}$.

One similarly observes that the eigenvalues of $\Delta_3^2$ for the $\K_3$ summands in  $V=\{\la_1^3\la_2^2\la_3\}$ coincide only if $\la_1^4=\la_2^2\la_3^2$, $\la_3^2=-\la_1\la_2$ or $\la_1^3=-\la_2^3$.
This again is not possible in $R/\p$ for $\p$ as in Table 2 for our given $V$. 

\section{Path representations}\label{paths:sec}

The main result of this section is the calculation of matrix entries for the braid generators and their eigenprojections for a path basis of a regular module. This will allow us to determine possible subquotient modules, in view of Lemma \ref{pbaseprop}.

\subsection{Affine braid group $AB_2$} 
The affine braid group $AB_2$ is defined via generators $\tau$ and $\sigma$ and relation
$\sigma\tau\sigma\tau=\tau\sigma\tau\sigma$. We will review how any path representation of braid groups
will give us
representations of the affine braid group $AB_2$ on a finite dimensional vector space $W$. 
More precisely, we will obtain matrices $A$ and $T$ which satisfy the relation 
\begin{equation}\label{affinebraid}
ATAT=TATA
\end{equation}
\begin{equation}\label{deltarelation}
ATAT=\delta 1_W,
\end{equation}
acting on a vector space $W$ with a basis $\{ v_r\}$ such that
\begin{equation}\label{Tdiagonal}
Tv_r=x_rv_r,\quad x_i\neq x_j\ {\rm for}\ i\neq j.
\end{equation}
It follows from \ref{deltarelation} and \ref{Tdiagonal} that 
\begin{equation}\label{Ainverse}
(A^{-1})_{rs}=\frac{x_rx_s}{\delta} a_{rs}.
\end{equation} 

\subsection{Generalized quantum Jucys-Murphy elements}
We recall how one can obtain representations of $AB_2$ from representations of ordinary braid groups using
(quantum) Jucys-Murphy elements. This was first done for Hecke algebras in \cite{Cherednik}.

\begin{lemma}\label{affinebraidlemma} Let $V$ be a $B_{n+1}$ module with a path basis, labeled by $\nu$.
Let $\la$ be a label for a simple $B_{n-1}$ module which appears if $V$ is viewed as a $B_{n-1}$ module,
and let $W(\la,\nu)$ be the vector space with a basis labeled by all paths of length 2 from $\la$ to $\nu$.
Then the assignment
$$\tau\mapsto T=(\Delta_n^2)_{|W(\la,\nu)},\hskip 3em \sigma\mapsto A=(\sigma_n)_{|W(\la,\nu)},$$
makes $W(\la,\nu)$ into an $AB_2$ module satisfying \ref{affinebraid},  \ref{deltarelation}, with $\delta=\Delta_{n-1}^2(\la)\Delta_{n+1}^2 (\nu)$,
and \ref{Tdiagonal}, except possibly that some of the $x_r$s coincide;  here $x_r$ is given by the scalar $\Delta_n^2(\mu_r)$ for the path $r: \la\to\mu_r\to\nu$,
and $\Delta_k^2(\ga)$ is the scalar via which $\Delta_k^2$ acts on a simple $B_k$ module labeled by $\ga$.
\end{lemma}

$Proof.$ This is well-known. We give some details to the readers who want to prove it for themselves.
It is probably easiest shown using pictures as in e.g. \cite{KT}
 that $\gamma_n=\Delta_n^2\Delta_{n-1}^{-2}=\sigma_{n-1}\ ...\ \sigma_2\sigma_1^2\sigma_2\ ...\ \sigma_{n-1}$.
We then have
\begin{equation}\label{qMurphy}
\sigma_n\gamma_n\sigma_n\gamma_n=\gamma_{n+1}\gamma_n=\gamma_n\gamma_{n+1}=\gamma_n\sigma_n\gamma_n\sigma_n.
\end{equation}
 Hence 
$$\tau\mapsto \gamma_n,\hskip 3em \sigma\mapsto \sigma_n,$$
defines a group homomorphism from $AB_2$ into $B_{n+1}$. As $\Delta_{n-1}^2$ commutes with both $\gamma_n$ and 
with $\sigma_n$, this remains a group homomorphism if we replace $\gamma_n$ by $\gamma_n\Delta_{n-1}^2=\Delta_n^2$.
To prove the claimed value of $\delta$ for \ref{deltarelation}, observe that $\gamma_{n+1}\gamma_n=\Delta_{n+1}^2\Delta_{n-1}^{-2}$, which acts via the scalar $\al\al'^{-1}$
on $W$. Hence $\sigma_n\Delta_n^2\sigma_n\Delta_n^2=\gamma_{n+1}\gamma_n\Delta_{n-1}^4$ acts via the scalar
$\al\al'^{-1}\al'^{2}$, as claimed.

\vskip .2cm



\subsection{Representations with $A$ having at most 2 eigenvalues}  
We consider representations of $AB_2$ on a vector space $W$ satisfying \ref{affinebraid}, \ref{deltarelation}, \ref{Tdiagonal} and $(A-\la_11)(A-\la_21)=0$.
The following is a version of the well-known quantum Jucys-Murphy construction, see \cite{Cherednik}

\begin{lemma}\label{2by2matrixlem} (a) If $A$ acts via the scalar  $\la$ on a basis vector $v_r$, then $\delta=x_r^2\la^2$. Hence we only have two possible values $x_r=\pm \sqrt{\delta}/\la$ of $T$ for given $\delta$.

(b) If $A$ only has two eigenvalues $\la_1$ and $\la_2$, with $\la_1$ having multiplicity 1, then $\dim W\leq 4$. We have at most one irreducible 2-dimensional subrepresentation span$\{ v_r,v_s\}$ of $AB_2$. 
The diagonal entries $d_r(\la_1)$ and $d_s(\la_1)$ of the  eigenprojection of $A$ for its
eigenvalue $\la_1$ with respect to this basis are given by
\begin{equation}\label{2dimformula}
d_r(\la_1)=- \frac{\la_1x_s+\la_2x_r}{(x_r-x_s)(\la_1-\la_2)}, \hskip 3em d_s(\la_1)=\frac{\la_1x_r+\la_2x_s}{(x_r-x_s)(\la_1-\la_2)}.
\end{equation}
\end{lemma}

$Proof.$ Part (a) is obvious. 
For part (b), it follows from \ref{Ainverse} and
$A+\la_1\la_2A^{-1}=(\la_1+\la_2)I$ that
$$(\delta + \la_1\la_2x_rx_s)a_{rs}=\delta(\la_1+\la_2)\delta_{rs}.$$
Hence $a_{rs}=0$ or $\delta + \la_1\la_2x_rx_s=0$ for $r\neq s$. As $x_s\neq x_{s'}$ for $s\neq s'$ by assumption,
there is at most one $s$ for given $r$ for which this is the case.
As $a_{rs}\neq 0$ for these $r,s$, we deduce from this  that
$$\delta = -\la_1\la_2x_rx_s.$$
and 
$$a_{rr}=\frac{\delta(\la_1+\la_2)}{\delta + \la_1\la_2x_r^2}, \hskip 3em a_{ss}=\frac{\delta(\la_1+\la_2)}{\delta + \la_1\la_2x_s^2}.$$
By definition, $d_r(\la_1)$ is the $r$-th diagonal entry of the eigenprojection $P=(A-\la_2I)/(\la_1-\la_2)$
of $A$ for the eigenvalue $\la_1$. We obtain the claimed result from this via a straightforward calculation.
The restriction of $A$ to span $\{ v_r, v_s\}$ has determinant $\la_1\la_2$ from which we calculate
$$a_{rs}a_{sr}=\frac{-(\la_2x_r+\la_1x_s)(\la_2x_s+\la_1x_r)}{(x_r-x_s)^2}.$$
Hence we obtain a 2-dimensional irreducible subrepresentation of $W$ spanned by $v_r$ and $v_s$. 
As $\la_1$ appears with multiplicity 1 in $W$, the complement of this subrepresentation is contained in the eigenspace of $A$ with eigenvalue $\la_2$. By (a), we only have two possibilities for the eigenvalue $x_r$ for the basis vector of a one-dimensional representation of $AB_2$.

\subsection{Examples}\label{n2examples} We will need the following examples of two-dimensional representations of $AB_2$. It will be convenient
in one case to call the second eigenvalue $\la_3$ in stead of $\la_2$, so we also list the eigenvalues of $A$.

\[
\begin{tabular}{|c|c|c|c|c|c|}paths& $x_1$ & $x_2$&e.values&$d_1(\la_1)$& $d_2(\la_1)$ \\
\hline
&& & & &\\
$\{1\}\to\{\la_1\la_2\}$&$\la_1^2$ & $\la_2^2$&$ \la_1,\la_2$& $\frac{-\la_1\la_2}{(\la_1-\la_2)^2}$& $\frac{\la_1^2-\la_1\la_2+\la_2^2}{(\la_1-\la_2)^2}$\\
$\{\la_2\}\to\{\la_1^3\la_2^2\la_3\}$&$(\la_1\la_2\la_3)^2$&$-(\la_1\la_2)^3$&
$ \la_1,\la_2$& $\frac{\la_2(\la_1^2-\la_3^2)}{(\la_1-\la_2)(\la_3^2+\la_1\la_2)}$&  $\frac{\la_1(\la_3^2-\la_2^2)}{(\la_1-\la_2)(\la_3^2+\la_1\la_2)}$\\
$\{\la_3\}\to\{\la_1^4\la_2^2\la_3^2\}$&$(\la_1\la_2\la_3)^2$&$-(\la_1\la_3)^3$&
$ \la_1,\la_2$& $\frac{\la_1^2\la_3-\la_2^3}{(\la_1-\la_2)(\la_2^2+\la_1\la_3)}$ & $ \frac{\la_1\la_2(\la_2-\la_3)}{(\la_1-\la_2)(\la_2^2+\la_1\la_3)}$ \\
$\{\la_2\}\to\{\la_1^4\la_2^2\la_3^2\}$&$(\la_1\la_2\la_3)^2$&$-(\la_1\la_2)^3$& 
$ \la_1,\la_3$&
$\frac{\la_1^2\la_2-\la_3^3}{(\la_1-\la_3)(\la_3^2+\la_1\la_2)}$& $\frac{\la_1\la_3(\la_3-\la_2)}{(\la_1-\la_3)(\la_3^2+\la_1\la_2)}$\\
\end{tabular}
\]
\vskip .2cm
Here $\{1\}\to \{\la_1\la_2\}$ stands for the two paths for the two-dimensional representation $\{\la_1\la_2\}$ of $K_3$, and $\ga\to\nu$ stands for the two paths of length 2
from the $K_2$ representation $\la$ to the $K_4$ representation $\nu$ in lines 2 to 4.

\subsection{More complicated representations of $AB_2$}\label{AB2section}
We now consider representations of $AB_2$ where $A$ satisfies a cubic equation, and, moreover, one of its eigenvalues, say $\la_3$,  has multiplicity 1.
So if $P$ is the corresponding eigenprojection for $\la_3$,  the matrices $A$ and $T$ satisfy, besides
relations \ref{affinebraid}, \ref{deltarelation} and \ref{Tdiagonal}, also the relation
\begin{equation}\label{cubicA}
(A-\lambda_1I)(A-\lambda_2I)(A-\lambda_3I)=0 \quad {\rm and\quad rank}\ P=1.
\end{equation}
 It follows from these conditions that
\begin{equation}\label{APrelation}
A+\lambda_1\lambda_2A^{-1}=(\la_1+\la_2)I + \la_3^{-1}(\la_1-\la_3)(\la_2-\la_3)P.
\end{equation}
We deduce from the previous equation and \ref{Ainverse} that
\begin{equation}\label{linequa}
\frac{\delta + \la_1\la_2x_ix_j}{\delta} a_{ij} = (\la_1+\la_2)\delta_{ij} + \la_3^{-1}(\la_1-\la_3)(\la_2-\la_3)p_{ij}.
\end{equation}

The following lemma and  theorem are slight improvements of \cite{Wexc}, Proposition 5.7, given in notation more suitable for our purposes.
\begin{lemma}\label{AB2lemma} Assume $W$ is an $AB_2$ module satisfying conditions \ref{affinebraid}, \ref{deltarelation} and \ref{Tdiagonal}, with $n\geq 2$. Moreover, we assume $A$ has at most three eigenvalues, one of which, say $\mu$, has multiplicity 1. If $\delta\neq \la_k^2x_r^2$ for all eigenvalues $\la_k\neq \mu$, and $1\leq r\leq n$, then $A$ and $T$ generate the full $n\times n$ matrix algebra.
\end{lemma}

$Proof.$ Let us first do the more complicated case with $A$ having three eigenvalues, with $\mu=\la_3$.
As $P$ has rank 1, we can write it as $P=vw^T$ for two vectors $v,w\in F^n$. Assume $v_i=0$. Then it follows from \ref{linequa} and our assumption $\delta+\la_1\la_2x_ix_j\neq 0$ that $a_{ij}=0$ for $j\neq i$. Hence $a_{ii}$ must be an eigenvalue of $A$. It can not be equal to $\la_3$ as the $i$-th column of its eigenprojection $P$ would be equal to 0. Hence $a_{ii}\in\{\la_1,\la_2\}$. If $a_{ii}=\la_1$, it follows from \ref{linequa} that $\la_1(1+\delta^{-1}\la_1\la_2x_i^2)=\la_1+\la_2$,
from which one deduces $\la_1^2x_i^2=\delta$. One similarly shows that $a_{ii}=\la_2$ implies $\la_2^2x_i^2=\delta$. One proves that $w_i=0$ implies $\la_k^2x_i^2=0$ for $k\in \{1,2\}$ the same way.
We can therefore assume that $v_rw_r\neq 0$ for $1\leq r\leq n$. Hence all entries of $P$ are nonzero.
Moreover, as all the $x_i$s are distinct, we can also obtain the diagonal matrix entries as eigenprojections of $T$. This proves that $A$ and $T$ generate the full $n\times n$ matrix algebra.

If $A$ only has two eigenvalues, the claim follows from Lemma \ref{2by2matrixlem}.

\vskip .3cm
\begin{theorem}\label{AB2theorem}(generalized $q$-Jucys-Murphy, see \cite{Wexc}) We assume a representation of $AB_2$ satisfying conditions  \ref{affinebraid}, \ref{deltarelation}, \ref{Tdiagonal} and \ref{cubicA} with $(\la_1-\la_3)(\la_2-\la_3)\neq 0$. 
If $\delta+\la_1\la_2x_rx_s\neq 0$ for $1\leq r,s\leq n$ and $\delta\not\in\{\la_k^2x_r^2,\ k=1,2,\ 1\leq r\leq n\}$, then  

(a) the diagonal entries $d_r(\la_3)$ of $P$ are nonzero and

(b) they are given by
$$d_ r(\la_3)\ =\ \frac{\la_3\ b(r,3)}{(\la_1-\la_3)(\la_2-\la_3)}\ \prod_{s\neq r}\frac{\delta+\la_1\la_2x_rx_s}{\delta(x_r-x_s)},$$
where
$$
b(r,3)\ =\ (-1)^{n-1}\la_3\frac{\la_1\la_2x_r+\delta x_r^{-1}}{\delta}\ (\prod_{t=1}^n x_t)
- (\la_1+\la_2)(\frac{-\delta}{\la_1\la_2})^{(n-1-\ve)/2}x_r^\ve,
$$
and $\ve=0$ for $n$ odd and $\ve=1$ for $n$ even. If all $d_r(\la_3)\neq 0$, our $AB_2$ representation is irreducible, and  $A=(a_{rt})$ is given uniquely, up to conjugation by a diagonal matrix by 
$$a_{rt}=\frac{b(r,3)}{x_r-x_t}\ \prod_{s\neq r,t}\frac{\delta+\la_1\la_2x_rx_s}{\delta(x_r-x_s)}.$$
\end{theorem}

$Proof.$ 
This is essentially just a version of the proof of \cite{Wexc}, Proposition 5.7 
in our more general
notation. We outline the main steps for the readers' convenience. Let $P=vw^T$ be as in the proof of Lemma \ref{AB2lemma}. As $v_i\neq 0$,
it follows from $Av=\la_3v$, after solving for $a_{ij}$ in Eq \ref{linequa}, that
$$\sum_j\ \frac{\delta}{\delta+\la_1\la_2x_ix_j}( (\la_1+\la_2)\delta_{ij}\ +\ \frac{(\la_1-\la_3)(\la_2-\la_3)}{\la_3}v_iw_j)v_j
\ =\ \la_3v_i.$$
Dividing by $v_i$ and setting $d_j=v_jw_j$, we obtain the linear system
$$\sum_j\ \frac{\delta(\la_3-\la_1)(\la_3-\la_2)}{\la_3(\delta + \la_1\la_2x_ix_j)}\ d_j\ =
\ \la_3-\frac{\delta(\la_1+\la_2)}{\delta + \la_1\la_2x_i^2}.$$
This can be transformed to
$$\sum_j\ \frac{1}{1+\la_1\la_2x_ix_j/\delta}\ d_j\ =
\ \frac{\la_3^2}{(\la_1-\la_3)(\la_2-\la_3)}\ -\ \frac{\la_3(\la_1+\la_2)}{(\la_1-\la_3)(\la_2-\la_3)(1 + \la_1\la_2x_i^2/\delta)}.$$
This system can now be solved after substituting $\tilde x_i = (-\la_1\la_2/\delta)^{1/2}x_i$, using 
Lemmas 5.5 and 5.6 from \cite{Wexc}. This yields the expressions for $d_r(\la_3)$ as stated. 

If all $d_r(\la_3)\neq 0$, the rank 1 matrix $P$ is uniquely determined by its diagonal entries, up to conjugation by a diagonal matrix.
Setting $P=(d_r(\la_3))_{rs}$, we obtain the matrix $A$ from \ref{linequa} as claimed.


\subsection{Examples}\label{examplesnthree} We now give some examples of $AB_2$ representations of dimension 3, which will be needed later. In this case, the formula for the $r$-th diagonal entry $d_r(\la_3)$ of the eigenprojection $P$ of $A$  for the eigenvalue $\la_3$ becomes
$$
d_ r(\la_3)\ =\ \frac{\la_3}{(\la_1-\la_3)(\la_2-\la_3)}\  [\la_3\frac{\la_1\la_2x_r+\delta x_r^{-1}}{\delta}\ (\prod_{t=1}^3 x_t)+
 (\la_1+\la_2)\frac{\delta}{\la_1\la_2})]\ \prod_{s\neq r}\frac{\delta+\la_1\la_2x_rx_s}{\delta(x_r-x_s)}.
$$

(1) We first consider the 3-dimensional representation of $\K_3$, with three paths  $\{1\}\to \{\la_1\la_2\la_3\}$ of length 2. We have eigenvalues $x_i=\la_i^2$ of $T$ for $1\leq i\leq 3$, and $\delta=(\la_1\la_2\la_3)^2$. We  then calculate the $i$-th diagonal entry 
$d_i(\la_j)$ of the eigenprojection of $A$ for the eigenvalue $\la_j$ as
$$d_i(\la_j)=\frac{(\la_i^2+\la_j\la_k)(\la_j^2+\la_i\la_k)}{(\la_i-\la_j)^2(\la_i-\la_k)(\la_k-\la_j)},\quad i\neq j,$$
where $\{ i,j,k\}=\{ 1,2,3\}$, and 
$$d_i(\la_i)=\frac{(\la_j+\la_k)^2\la_i^2}{(\la_j-\la_i)^2(\la_i-\la_k)^2}.$$
This can be checked explicitly for $d_i(\la_3)$ using Theorem \ref{AB2theorem}. The other cases are obtained from this using the permutation automorphisms in  Lemma \ref{permautlem}.

(2) We next consider the case for the three paths $\{\la_1\}\to\{\la_1^3\la_2^2\la_3\}$, with $x_1=(\la_1\la_2\la_3)^2$, $x_2=-(\la_1\la_2)^3$, $x_3=\la_1^6$ and
$\delta = \la_1^8\la_2^4\la_3^2$.
Let $d_r(\la_i)$ be the $r$-th diagonal entry of the eigenprojection of $A$ for the eigenvalue $\la_i$.
Moreover, our indices satisfy $\{ i,j,k\} = \{ r,s,t\} = \{ 1,2,3\}$. Then the formula 
for the diagonal entries simplifies to 
$$
d_r(\la_ i)\ =\ \frac{\la_i\la_1^4\la_2^2\la_3}{(\la_i-\la_j)(\la_i-\la_k)}
\  [ (\la_j^{-1}+\la_k^{-1})\la_1^4\la_2^2\la_3-x_r-\la_1^7\la_2^3\la_3\la_i x_r^{-1})]
\ \prod_{s\neq r}\frac{1+\delta^{-1}\la_j\la_kx_rx_s}{x_s-x_r}.
$$
We calculate the expression in straight brackets above
 for all values of $r$ and $i$ below. E.g. for $i=2$ the formula above would
simplify to $(\la_1+\la_3)\la_1^3\la_2^2-x_r-\la_1^7\la_2^4\la_3x_r^{-1}$.
\[
\begin{tabular}{|c|c|c|c|} $r\backslash i$ & 1&2&3 \\
\hline
&&&\\
1 & $\la_1^2\la_2\la_3^{-1}(\la_3^2-\la_1^2)(\la_1^2-\la_2\la_3)$&
$ -\la_1^2\la_2^2\la_3^{-1}(\la_1-\la_3)^2(\la_1+\la_3)$&$-\la_1^2\la_2(\la_1^2-\la_2\la_3)(\la_1-\la_3)$\\
2&$\la_1^3(\la_2^2+\la_1\la_3)(\la_1+\la_2)$&
$\la_1^3\la_2(\la_1+\la_2)(\la_2+\la_3)$&$\la_1^3(\la_2^2+\la_1\la_3)(\la_2+\la_3)$\\
3&$-\la_1^2(\la_1^2-\la_2^2)(\la_1^2-\la_2\la_3)$&
$-\la_1(\la_1^2-\la_2^2)(\la_1^3-\la_2^2\la_3)$&$-\la_1(\la_1^3-\la_2^2\la_3)(\la_1^2-\la_2\la_3)$\\
\end{tabular}
\]
We also have
$$x_1-x_2=\la_1^2\la_2^2(\la_3^2+\la_1\la_2),\quad x_2-x_3=-\la_1^3(\la_2^3+\la_1^3),\quad x_3-x_1=\la_1^2(\la_1^4-\la_2^2\la_3^2).$$
In the following table we calculate the quantity $1+\la_j\la_kx_rx_s/\delta= 1-\la_1^4\la_2^2\la_3/\la_ix_t$, where $\{ r,s,t\}=\{ 1,2,3\}$.

\[
\begin{tabular}{|c|c|c|c|}$i\backslash t$ & 1 &2&3\\ 
\hline
1 &$1-\la_1\la_3^{-1}$ &$1+\la_3\la_2^{-1}$&$1-\la_1^{-3}\la_2^2\la_3$\\ 
2 & $1-\la_1^2\la_2^{-1}\la_3^{-1}$ &$1+\la_1\la_2^{-2}\la_3$&$1-\la_1^{-2}\la_2\la_3$\\
 3 & $1-\la_1^2\la_3^{-2}$ &$1+\la_1\la_2^{-1}$&$1-\la_1^{-2}\la_2^2$\\
\end{tabular}
\]
For $i=1$ we obtain the values 
\begin{align}
d_1(\la_1)\ &=\ \frac{(\la_1+\la_3)(\la_2+\la_3)(\la_1^3-\la_2^2\la_3)}{(\la_1-\la_2)(\la_1^2+\la_2\la_3)(\la_3^2+\la_1\la_2)}\cr
d_2(\la_1)\ & =\ \frac{-(\la_2^2+\la_1\la_3)(\la_1^3-\la_2^2\la_3)}{(\la_1-\la_2)(\la_3^2+\la_1\la_2)
(\la_1^2-\la_1\la_2+\la_2^2)},\cr
d_3(\la_1)\ &=\ \frac{\la_1^2\la_2(\la_2+\la_3)}{(\la_1^2+\la_2\la_3)(\la_1^2-\la_1\la_2+\la_2^2)}.
\end{align}
For $i=2$ we obtain the values 
\begin{align}
d_1(\la_2)\ &=\ \frac{-\la_2(\la_1+\la_3)(\la_1-\la_3)^2(\la_2^2+\la_1\la_3)}{(\la_1-\la_2)(\la_2-\la_3)(\la_1^2+\la_2\la_3)(\la_3^2+\la_1\la_2)},\cr
d_2(\la_2)\ &=\  \frac{\la_2(\la_2+\la_3)(\la_1^2-\la_2\la_3)^2}{(\la_1-\la_2)(\la_2-\la_3)(\la_3^2+\la_1\la_2)
(\la_1^2-\la_1\la_2+\la_2^2)},\cr
d_3(\la_2)\ &=\ \frac{-(\la_1^3-\la_2^2\la_3)(\la_2^2+\la_1\la_3)}{(\la_2-\la_3)(\la_1^2+\la_2\la_3)(\la_1^2-\la_1\la_2+\la_2^2)}.
\end{align}
For $i=3$ we obtain the values
\begin{align}
d_1(\la_3)\ &=\ \frac{\la_3^2(\la_1+\la_2)^2(\la_1-\la_2)}{(\la_2-\la_3)(\la_1^2+\la_2\la_3)(\la_3^2+\la_1\la_2)},\cr
d_2(\la_3) \ &=\  \frac{-(\la_2+\la_3)(\la_1-\la_2)(\la_1+\la_3)(\la_2^2+\la_1\la_3)}{(\la_3^2+\la_1\la_2)(\la_2-\la_3)
(\la_1^2-\la_1\la_2+\la_2^2)},\cr
d_3(\la_3)\ &=\ \frac{\la_2(\la_1+\la_3)(\la_1^3-\la_2^2\la_3)}{(\la_2-\la_3)(\la_1^2+\la_2\la_3)(\la_1^2-\la_1\la_2+\la_2^2)}.
\end{align}
(3) We  consider the case with three paths $\{\la_1\}\to\{\la_1^3\la_2^3\la_3^3\}_\theta$, where we have
$x_1=-(\la_1\la_2)^3,\quad x_2=-(\la_1\la_3)^3,\quad x_3=(\la_1\la_2\la_3)^2$
and $\delta=\theta\la_1^6\la_2^4\la_3^4$. We then obtain the diagonal entries (for the eigenprojection
of the eigenvalue $\la_3$)
\begin{align}
 d_1(\la_3)\ &=\ \frac{\theta(\la_1+\theta\la_3)(\la_2+\theta\la_3)(\la_3+\theta\la_2)(\la_3^2-\theta\la_1\la_2)}
{(\la_2-\la_3)^2(\la_3-\la_1)(\la_3^2+\la_1\la_2)},\cr
d_2(\la_3)\ &=\ \frac{\theta^{-1}\la_2\la_3(\la_1+\theta\la_3)(\la_2^2-\theta\la_1\la_3)}
{(\la_2-\la_3)^2(\la_1-\la_3)(\la_2^2+\la_1\la_3)},\cr
d_3(\la_3)\ &=\ \frac{\theta\la_1\la_3(\la_1+\theta\la_2)(\la_3-\theta\la_2)^2(\la_3+\theta\la_2)}
{(\la_1-\la_3)(\la_2-\la_3)(\la_2^2+\la_1\la_3)(\la_3^2+\la_1\la_2)}.
\end{align}
The diagonal entries for the eigenprojection of $\la_2$ are obtained from the ones just calculated 
by permuting $\la_2$ with $\la_3$ and then by permuting $d_1$ with $d_2$. The diagonal entries
for the eigenprojection of $\la_1$ are given by
\begin{align}
d_1(\la_1)\ &=\ \frac{\theta^{-1}\la_1\la_3(\la_2^2-\theta\la_1\la_3)(\la_2+\theta\la_3)}
{(\la_2-\la_3)(\la_3^2+\la_1\la_2)(\la_2-\la_1)(\la_1-\la_3)},\cr
d_2(\la_1)\ &=\ \frac{\theta^{-1}\la_1\la_2(\la_3^2-\theta\la_1\la_2)(\la_3+\theta\la_2)}
{(\la_2-\la_3)(\la_2^2+\la_1\la_3)(\la_2-\la_1)(\la_1-\la_3)},\cr
d_3(\la_1)\ &=\ \frac{-\theta(\la_2^2-\theta\la_1\la_3)(\la_3^2-\theta\la_1\la_2)(\la_1+\theta\la_2)(\la_1+\theta\la_3)}
{(\la_3^2+\la_1\la_2)(\la_2^2+\la_1\la_3)(\la_2-\la_1)(\la_1-\la_3)}.
\end{align}

\subsection{Example}\label{n4example} We consider the case with four paths  $\{\la_1\}\to\{\la_1^4\la_2^2\la_3^2\}$
of length 2. We then have $x_1=\la_1^6$, $x_2=-(\la_1\la_2)^3$, 
$x_3=(\la_1\la_2\la_3)^2$, $x_4=-(\la_1\la_3)^3$ and $\delta= \la_1^8\la_2^3\la_3^3$.
Plugging these values into the formula for $d_i$ in Theorem \ref{AB2theorem},
we obtain
$$d_r(\la_3)=\frac{x_r-\la_1^6\la_2^3\la_3^3x_r^{-1}}{\la_1^{13}\la_2^4\la_3^5(\la_3-\la_1)(\la_3-\la_2)}
\prod_{s\neq r}\frac{\la_1^7\la_2^2\la_3^3+x_rx_s}{(x_r-x_s)}.$$
Plugging the values of $x_i$ into these formulas, we obtain
\begin{align}\label{n4sol}
d_1(\la_3)\ &=\ \frac{(\la_1^2\la_2-\la_3^3)(\la_1^4+\la_1^2\la_2\la_3+\la_2^2\la_3^2)(\la_1-\la_2)}
{(\la_1-\la_3)(\la_2-\la_3)(\la_1^2+\la_2\la_3)(\la_1^2-\la_1\la_2+\la_2^2)(\la_1^2-\la_1\la_3+\la_3^2)},\cr
d_2(\la_3)\ &=\ \frac{-(\la_2^3-\la_1^2\la_3)(\la_3^3-\la_1^2\la_2)}
{(\la_1\la_2+\la_3^2)(\la_1-\la_3)(\la_2-\la_3)(\la_1^2-\la_1\la_2+\la_2^2)},\cr
d_3(\la_3)\ &=\ \frac{-\la_2\la_3(\la_1+\la_3)^2(\la_2^3-\la_1^2\la_3)}
{(\la_2-\la_3)(\la_1^2+\la_2\la_3)(\la_2^2+\la_1\la_3)(\la_3^2+\la_1\la_2)},\cr
d_4(\la_3)\ &=\ \frac{-\la_3^2(\la_1+\la_2)^2(\la_1-\la_2)}
{(\la_2^2+\la_1\la_3)(\la_1^2-\la_1\la_3+\la_3^2)(\la_2-\la_3)}.
\end{align}
\subsection{Determining matrix blocks for the braid generators}\label{matrixblocks}
We will encounter non-semisimple representations of $AB_2$ in the following context. Assume that $\delta$ and the eigenvalues $x_r$ of $T$ are monomials in the eigenvalues $\la_i$. Let $R=F_0[\la_1,\la_2,\la_3]$ and let $\p$ be a prime ideal in $R$. We assume as before $W$ to be an $AB_2$ module with a basis of eigenvectors of $T$ with mutually distinct eigenvalues. Moreover, we assume that the matrix coefficients of $A$ are in the localization $R_\p$. This implies that we can also get a representation over $Q(R/\p)$.

\begin{corollary}\label{nonsemisimp}
Assume that the diagonal entries $d_r$ of $P$ are not in $\p R_\p$ except for one entry $d_t$, which is not in $\p^2R_\p$.
Also assume that  the entries $\delta+\la_1\la_2x_rx_s\not\in \p$ for all indices $r,s$. Then there exist two non-equivalent non-semisimple representations of $AB_2$ over the quotient field $Q(R/\p)$, in which the basis vector $v_t$ either spans a 1-dimensional  sub module or quotient module.
\end{corollary}

$Proof.$ By assumption, we have $d_td_s\in \p R_\p\backslash \p^2R_\p$ for all $s\neq t$.
As $P$ is a rank 1 matrix, this implies $p_{ts}p_{st}=d_td_s\in \p R_\p\backslash \p^2R_\p$ for all $s\neq t$.
Hence, for given $s$, exactly one of the two entries  $p_{ts}$ or $p_{st}$ is in $\p R_\p$, say $p_{st}$.
This implies that $\bar p_{st}=0$ and $\bar p_{ts}\neq 0$ in $Q(R/\p)$.
But as $P$ has rank 1, we have $\bar p_{\tilde{s}t}=0$ for all $\tilde{s}\neq t$. It follows that $v_t$ spans a 1-dimensional $AB_2$ submodule
in the representation over $Q(R/\p)$. If $p_{ts}\in\p R_\p$, one similarly shows that $v_t$ spans a 1-dimensional quotient module.
\begin{remark}\label{n4remark} Assume we have eigenvalues $\la_i$ such that the matrix entries $d_r(\la_3)$, $1\leq r\leq 4$  above are well-defined, and define the projection $P$ via the matrix $(d_r(\la_3))_{rs}$. Then it is  straightforward, if moderately tedious to determine the entries of the braid matrix $A$ 
using the last formula in Theorem \ref{AB2theorem}.  One can thus check that the entry $a_{rs}$  is zero mod 
$(\la_2^3-\la_1^2\la_3)$
if and only if $r\in \{2,3\}$ and $s\in \{1,4\}$.  One similarly checks that the entries $a_{rs}$ are zero mod $\la_1^4+\la_1^2\la_2\la_3+\la_2^2\la_3^2$ if and only if $r=1$ and $s>1$.
\end{remark}

\section{Structure results for  $\K_4$}\label{structure:sec}

\subsection{Strategy}\label{strategy:sec}  Let $V$ be a regular module. We already know that it has to be simple unless its eigenvalues satisfy one of the polynomials listed in Proposition \ref{semisimplevalues:prop}.
We fix such a polynomial and the corresponding prime ideal $\p$. 
It follows from Lemma \ref{pathexistence} that our path representations are well-defined, if we consider $V$ as a $Q(R/\p)$ module for $\p$ one of the corresponding prime ideals; here $Q(R/\p)$ is the quotient field of $R/\p$. We now describe how knowledge of diagonal entries of eigenprojections of the braid generators in these representations will be used to determine the structure of these modules.
It may perhaps be more instructive following one of the examples worked out in the proofs of Lemmas \ref{8dimcomplem} or \ref{9dimcomplem}.

(1)  {\it Bases for simple subquotients} Observe that we have at least two rank 1 eigenprojections in each block of the matrix $S_i$ for a braid generator $\sigma_i$ in the path representation. We determine the equivalence classes of the relation
defined as the transitive closure of $s\leftrightarrow t$ if $d_sd_t\not\in \p$ for the diagonal entries  of at least one of these eigenprojections, see Lemma \ref{pbaseprop}.  
As we shall see, there will be at most two such equivalence classes. If we only have one such equivalence class, $V$ will be simple as a $Q(R/\p)$ module, see Lemma \ref{matrixcoeff}.

(2) {\it Composition series} If we have two equivalence classes, we are going to show that they correspond to a sub and a quotient module of $V$ as follows:

(a) Pick a matrix block for a braid generator (see Lemma \ref{pbaseprop}) of maximum size with paths in both equivalence classes. Calculate the braid matrix using Theorem \ref{AB2theorem},  Corollary \ref{nonsemisimp} or Remark \ref{n4remark}.
This depends on a choice about the off-diagonal entries of the eigenprojection $P$. We can choose the matrices $(d_r(\la_i))_{rs}$ or $(d_s(\la_i))_{rs}$. Replacing one choice by the other one will result in
interchanging sub and quotient modules of $V$. In the language of  Remark \ref{pathunique}, this would correspond to interchanging the module $V$ by $V^t$.

(b) If there are other matrix blocks with non-equivalent paths, put question marks in the corresponding off-diagonal entries.

(c) Use braid relations to determine which of the question mark off-diagonal entries have to be zero mod $\p$.

\ni This determines the structure of $V$ as a $Q(R/\p)$ module. We will deal with the structure of $V$ as a $Q(R/\q)$ module for ideals $\q$ larger than the ones treated so far in Section \ref{nomore:sec}.





\subsection{Subquotients for the 6-, 8- and 9-dimensional representations of $B_4$}\label{6dim2} We are going to determine the structure of these regular modules, viewed as $Q(R/\p)$ modules, as discussed in the previous subsection.

\begin{lemma}\label{6dimlem}
Let $\p$ be one of the ideals listed in Table 2 for the module $V=\{\la_1^3\la_2^2\la_3\}$. Then the $Q(R/\p)$ module $V$ or $V^t$ has simple quotient and submodules as listed in Corollary \ref{6dimquot}.
\end{lemma}

$Proof.$  The existence of the claimed sub and quotient modules  has already been shown in Corollary \ref{6dimquot}.
To prove simplicity of these modules, we can apply Lemma \ref{matrixcoeff}, either by using the matrix coefficients in 
Lemma \ref{6dimweightbasis}, or the diagonal entries of $P_3(\la_3)$ in Example (2) in Section \ref{examplesnthree}, corresponding to the paths
from $\{\la_1\}$ to $\{\la_1^3\la_2^2\la_3\}$. Also observe that if we have a nontrivial quotient module, there exists at least
one matrix block for $S_i$, $i\in\{ 2,3\}$ where some basis  paths belong to the quotient, and some to the submodule. This is possible only if one of the diagonal entries $d_t$ for a rank 1 eigenprojection of $S_i$ is in $\p R_\p$. One can check from the results in Section \ref{examplesnthree}, (2), that these diagonal entries are not in $\p^2R_\p$. But then the corresponding $AB_2$ representation does not split, see Corollary \ref{nonsemisimp}. Hence the $\K_3$ quotient can not split either.

\medskip

We will have to work harder to prove nonsemisimplicity of the 8- and 9-dimensional regular modules, as we are not aware of a similar
connection to quantum groups as in the previous lemma.

\begin{lemma}\label{8dimcomplem} Let $V=\{\la_1^4\la_2^2\la_3^2\}$.

(1) If $\p=(\la_2+\la_3)$ or $\p=(\la_1^2-\la_2\la_3)$, the $Q(R/\p)$ module $V$ is simple.


(2) If $\p=(\la_2^3-\la_1^2\la_3)$, the $Q(R/\p)$ module $V$ has a 5-dimensional quotient  $\{\la_1^2|\la_2^2\la_3\}$ with
weights $(11),(12),(21),$ $(23)$ and $(32)$, and a submodule isomorphic to the generic representation $\{\la_1^2\la_3\}$.
The case  $\p=(\la_3^3-\la_1^2\la_2)$ is similar.  

(3) If $\p=(\la_1^2-\theta\la_2\la_3)$, the $Q(R/\p)$ module $V$ has
 a 7-dimensional quotient representation $\{\la_1^3\la_2^2\la_3^2\}$  with the same weights as  $\{\la_1^4\la_2^2\la_3^2\}$, but with the weight $(1,1)$ only appearing once.
\end{lemma}

$Proof.$ Observe that $\la_1^6-\la_2^3\la_3^3=(\la_1^2-\la_2\la_3)(\la_1^2-\theta\la_2\la_3)(\la_1^2-\theta^{-1}\la_2\la_3)$. 
Hence, by  Lemma \ref{9and8dim}, the representation  $\{\la_1^4\la_2^2\la_3^2\}$ is simple except possibly if $\la_2=-\la_3$, $\la_1^2=\la_2\la_3$ or one of the cases in (2) or (3) here apply. Let $\p$ be the prime ideal generated by one of these irreducible
polynomials.
It follows from Table 1 and Theorem \ref{B3class} that in all of these cases we have the isomorphism of $Q(R/\p)\K_3$ modules 
\begin{equation}\label{decomp1}
\{\la_1^4\la_2^2\la_3^2\}\ \cong\ \{\la_1\} \oplus \{\la_1\la_2\}\oplus \{\la_1\la_2\la_3\}\oplus \{\la_1\la_3\},
\end{equation}
where each summand on the right hand side is simple.  Now observe that all diagonal entries $d_r(\la_3)$  of the eigenprojection $P$ in the $4\times 4$ block of $S_3$ are not in $\p$ if $\p=(\la_2+\la_3)$ or $\p=(\la_1^2-\la_2\la_3)$, see Section \ref{n4example}.
It follows that all paths for our module are equivalent, i.e. the $Q(R/\p)$-module
$\{\la_1^4\la_2^2\la_3^2\}$ is simple in these cases. This proves (1).

Let now $\p=(\la_2^3-\la_1^2\la_3)$.
We order the paths for $\{\la_1^4\la_2^2\la_3^2\}$ by
starting with the left most $\K_3$ module in \ref{decomp1}, and ordering the paths in a $\K_3$-module
with the ordering of the eigenvalues $\la_i$.
E.g the path corresponding to eigenvalue $\la_2$ for $\sigma_1$ in the $\K_3$ module $\{\la_1\la_2\la_3\}$ would
be the fifth path. See also the paths into $\La_1+\La_2$ in the figure below \cite{MWa}, Remark 2.2.

As $\K_3$ is semisimple over $Q(R/\p)$ in this case,
paths belonging to the same $\K_3$ module are equivalent. 
 The 4 by 4 block for the matrix $S_3$ for $\sigma_3$ can be calculated from the example  in Section \ref{n4example}.
We will use the notation in that section, where the base vectors $v_i$, $1\leq i\leq 4$ correspond to
the paths 1, 2, 4 and 7 here. We indicate by 0 those entries of $S_4$ which are 0 mod $\p R_\p$, see Remark \ref{n4remark}.
We similarly calculate the 2 by 2 blocks for $S_3$ for the paths  from $\{\la_2\}$ resp $\{\la_3\}$ to  $\{\la_1^4\la_2^2\la_3^2\}$, 
using the results in Section \ref{n2examples}. Again, one of the diagonal entries for the path from $\{\la_3\}$ is zero  mod $\p R_\p$.
This implies that one of the two entries with $?$ has to be equal to 0  mod $\p R_\p$, see Corollary \ref{nonsemisimp}.
$$
S_2\ =\ \left[\begin{matrix} x&&&&&&&\cr &x&x&&&&&\cr &x&x&&&&&\cr &&&x&x&x&& \cr &&&x&x&x&& \cr &&&x&x&x&& \cr
&&&&&&x&x\cr &&&&&&x&x\cr \end{matrix}\right]\hskip 3em
S_3\ =\ \left[\begin{matrix} x&x&&x&&&x&\cr 0&x&&x&&&0&\cr &&x&&x&&&\cr 0&x&&x&&&0& \cr &&x&&x&&& \cr &&&&&x&&? \cr
 x&x&&x&&&x&\cr &&&&&?&&x\cr \end{matrix}\right]
$$
Now observe that $S_2S_3S_2v_1\in {\rm span}\{ v_1, v_7, v_8\}$. But if the (6,8) entry $a_{68}^{(3)}$ of $S_3$ was nonzero,
$S_3S_2S_3v_1$ would also have a nonzero coefficient for $v_6$. Hence  $a_{68}^{(3)}=0$  mod $\p R_\p$. 
It follows that the image of $\K_4$ has a submodule spanned by $\{ v_1, v_7,v_8\}$ and a quotient spanned by
the remaining basis vectors.

Part (3) with  $\p=(\la_1^2-\theta\la_2\la_3)$ can be shown similarly, and is much easier. Again, calculating the matrix $S_3$ as in the previous case,
it follows from Remark \ref{n4remark} that all non-diagonal entries of $S_3$ in the first row are equal to 0  mod $\p R_\p$.
Hence the vectors $v_2$ to $v_7$ span a submodule.
This submodule is irreducible, as it is semisimple as a $\K_3$ module, and all other diagonal entries of $P_3(\la_3)$
in the 4 by 4 block are nonzero   mod $\p R_\p$.

\begin{lemma}\label{9dimcomplem}  We only get  nontrivial quotients of the 9-dimensional representation $\{\la_1^3\la_2^3\la_3^3\}_\theta$ over $Q(R/\p)$ for a prime ideal $\p$ in the following cases,
where $\{ i,j,k\}=\{ 1,2,3\}$:

(1) If  $\p=(\la_i^2-\theta\la_j\la_k)$, we get a 7-dimensional quotient with the same weights as the subquotient appearing in the 8-dimensional representation.

(2) If $\p=(\la_i+\theta\la_j)$, we have 6- and 3-dimensional subquotients with the same weights as the generic modules $\{\la_j^3\la_k^2\la_i\}$ and $\{\la_i^2\la_k\}$.
\end{lemma}

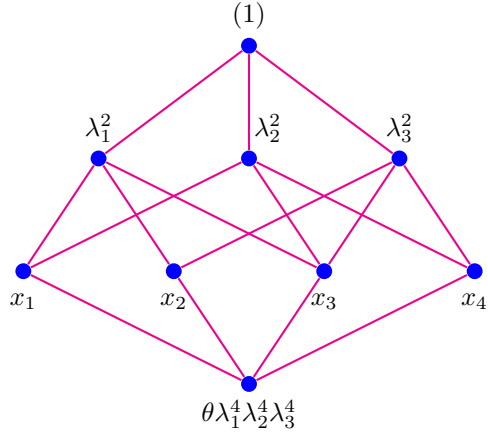
\begin{figure}[h]
\centering
\begin{tikzpicture}[
  dot/.style={circle, fill=blue, minimum size=6pt, inner sep=0pt},
  font=\small
]

\node[dot] (top) at (6,6) {};
\node at (6,6.4) {$(1)$};

\node[dot] (l1) at (4,4.5) {};
\node at (4,4.9) {$\lambda_1^2$};

\node[dot] (l2) at (6,4.5) {};
\node at (6.25,4.9) {$\lambda_2^2$}; 

\node[dot] (l3) at (8,4.5) {};
\node at (8,4.9) {$\lambda_3^2$};

\node[dot] (x1) at (3,3) {};
\node at (3,2.6) {$x_1$};

\node[dot] (x2) at (5,3) {};
\node at (5,2.6) {$x_2$};

\node[dot] (x3) at (7,3) {};
\node at (7,2.6) {$x_3$};

\node[dot] (x4) at (9,3) {};
\node at (9,2.6) {$x_4$};

\node[dot] (bottom) at (6,1.5) {};
\node at (6,1.1) {$\theta\lambda_1^4 \lambda_2^4\lambda_3^4$};

\draw[magenta, thick] (top) -- (l1);
\draw[magenta, thick] (top) -- (l2);
\draw[magenta, thick] (top) -- (l3);

\draw[magenta, thick] (bottom) -- (x1);
\draw[magenta, thick] (bottom) -- (x2);
\draw[magenta, thick] (bottom) -- (x3);
\draw[magenta, thick] (bottom) -- (x4);

\draw[magenta, thick] (l1) -- (x1);
\draw[magenta, thick] (l1) -- (x2);
\draw[magenta, thick] (l1) -- (x3);

\draw[magenta, thick] (l2) -- (x1);
\draw[magenta, thick] (l2) -- (x3);
\draw[magenta, thick] (l2) -- (x4);

\draw[magenta, thick] (l3) -- (x2);
\draw[magenta, thick] (l3) -- (x3);
\draw[magenta, thick] (l3) -- (x4);

\end{tikzpicture}
\caption{Paths for the 9-dimensional representation $\{\la_1^3\la_2^3\la_3^3\}_\theta$.}
\label{fig:bratteli-9d}
\end{figure}

$Proof.$ We have shown in Lemma \ref{9and8dim} that the 9-dimensional representation $V=\{\la_1^3\la_2^3\la_3^3\}_\theta$
is simple, except possibly if $\la_i=-\theta\la_j$ or $\la_i^2=\theta\la_j\la_k$. 
We will again use path representations to study these cases, as indicated in figure \ref{fig:bratteli-9d}. Here a vertex at level $k$ indicates a representatioon of $\K_k$, with the scalar via which $\Delta_k^2$ acts on it.  So we have $x_1 = -(\lambda_1 \lambda_2)^3$, $x_2 = -(\lambda_1 \lambda_3)^3$, $x_3 = (\lambda_1 \lambda_2 \lambda_3)^2$, and $x_4 = -(\lambda_2 \lambda_3)^3$.

Let us consider the case $\p=(\la_1+\theta\la_2)$. It follows from Section \ref{examplesnthree}, Example (3) and  Lemma \ref{permautlem} that all matrix entries of the eigenprojections
$P_3(\la_i)$ are well-defined. Hence we can use the same method as in Lemma \ref{8dimcomplem} to determine subquotient modules.
We label the paths as in  Lemma \ref{8dimcomplem} based on the isomorphism of $\K_3$ modules
$$\{\la_1^3\la_2^3\la_3^3\}_\theta\ \cong\ \{\la_1\la_2\}\oplus \{\la_1\la_3\}\oplus \{\la_1\la_2\la_3\}\oplus \{\la_2\la_3\}.$$

This determines the block structure for the matrix $S_2$. The matrix $S_3$ consists of three 3 by 3 blocks, which live on span $\{ v_1, v_3,v_5\}$,
span $\{ v_2, v_6,v_8\}$ and span $\{ v_4, v_7,v_9\}$.  The diagonal entries of the eigenprojections $P_3(\la_3)$ and $P_3(\la_1)$ for the first block are calculated in
 Section \ref{examplesnthree}, (3), and the ones for $P_3(\la_2)$ are obtained by applying the permutation $(23)$.  It follows that the diagonal
entries of $P_3(\la_1)$ and $P_3(\la_3)$ for the vector $v_5$ are equal to zero   mod $\p R_\p$.

One similarly gets the diagonal entries of the eigenprojections for the block acting on span $\{ v_4, v_7,v_9\}$ by applying the permutation $(13)$, see Lemma \ref{permautlem}.
One obtains that the diagonal entries of the projections $P_3(\la_2)$ and $P_3(\la_3)$ for the vector $v_4$ are equal to 0. Finally, we obtain from
Section  \ref{n2examples} that at least one of the two off-diagonal entries for the first 2 by 2 block of $S_2$ has to be equal to 0.

We now proceed as in the proof of Lemma \ref{8dimcomplem} as follows: We can set the non-diagonal entries on the 5-th column of $S_3$ equal to 0.
We then put question marks into the nondiagonal entries in the fourth row and fourth column of $S_3$ in the $(4,7,8)$ block and
into the non-diagonal entries for the $(12)$ block of $S_2$.  It is then straightforward to check that $S_3S_2S_3v_5$ is in the span of all basis vectors except $v_1$,
while this is true for $S_2S_3S_2v_5$ only if $(S_2)_{12}=0$  mod $\p R_\p$. 
If the non-diagonal entries of the fourth column of $S_3$ were 0  mod $\p R_\p$, then
$S_3S_2S_3v_4$ would be in the span of $\{ v_1, v_3, v_4, v_5\}$, while $S_2S_3S_2v_4$ is in the span of the basis vectors
$v_i$ with $i\in \{ 2,5,6,7,8,9\}$. As this would contradict the braid relation, the non-diagonal entries of $S_3$ in its fourth row must be equal to 0  mod $\p R_\p$. This proves that we get a $\K_4$ submodule
spanned by the vectors $v_i$ with $i\in \{ 2,5,6,7,8,9\}$ and a quotient module spanned by the remaining basis vectors.
As  $\la_1^2+\la_3^2\neq 0$  mod $\p R_\p$, the 3-dimensional subquotient is isomorphic to $\{\la_1^2\la_3\}$, which has the weights $(11),(13)$ and $(31)$. This also determines the weights of the 6-dimensional subquotient.
The general case $\p=(\la_i+\theta\la_j)$ follows from Lemma  \ref{permautlem} by symmetry.

 If  $\p=(\la_3^2-\theta\la_1\la_2)$, we can assume $\{\la_1^3\la_2^3\la_3^3\}_\theta$ to be a semisimple $\K_3$ module.
We leave it to the reader to check that one can set the nondiagonal entries of $S_3$ in the first row equal to 0, and that this
will also force that its nondiagonal entries in the second row will have to be equal to zero  mod $\p R_\p$.
This implies that we have a quotient module spanned by the vectors $v_1$ and $v_2$, and a submodule spanned by the
remaining basis vectors.

\subsection{Nonsemisimplicity of $\K_4$} We can now prove the converse of Proposition \ref{semisimplevalues:prop}.

\begin{theorem}\label{semisimple:thm} The regular module $V$ is not semisimple if and only if the eigenvalues of $\sigma_1$ are in the locus of one of the ideals listed next to it in Table 2, or if two of the eigenvalues coincide.
In particular, the algebra $\K_4$ is not semisimple if and only if  either $\sigma_1$ is not diagonalizable, or if its eigenvalues are in the locus of one of the ideals  listed  in Table 2.
\end{theorem}

$Proof.$ It remains to prove that the algebras $\K_4$ are not semisimple if one of the polynomials in  Table 2 is zero at the eigenvalues of $\sigma_1$. For this, it is enough to find a representation which is not semisimple. This is trivially true if two of the eigenvalues coincide, by mapping all $\sigma_i$ to the corresponding Jordan form.
For the polynomials $\la_i+\theta\la_j$ and $\la_i\pm \sqrt{-1}\la_j$, the result is known from the study of Iwahori Hecke algebras. For the other polynomials, we have shown in Lemmas \ref{6dimlem},  \ref{8dimcomplem}  and   \ref{9dimcomplem}
that our regular modules have nontrivial quotients which do not split if the eigenvalues satisfy one of the polynomials listed in Table 2. 

\begin{lemma}\label{uniqueness} Assume that the eigenvalues of $\sigma_i$ satisfy at most one of the polynomials in Table 2. If two simple $\K_4$ representations are isomorphic as $\K_3$ modules, and they have the same weights, then they are isomorphic, except for the two non-isomorphic 9-dimensional modules.
\end{lemma}
$Proof.$ This follows for the generic modules by inspection of Table 1. However, it might be conceivable that e.g. the 7-dimensional modules appearing as subquotients of the 8- and 9-dimensional generic representations might not be isomorphic. To rule that out,
observe that they have isomorphic path bases; for the subquotient of the nine-dimensional representation one obtains it by removing  the paths going through $\{\la_2\la_3\}$ in the picture in the proof of 
Lemma \ref{9dimcomplem}. We will denote the path going through the representation $\{\la_i\}$ and the $\K_3$ representation on which $\Delta_3^2$ acts via $x_r$ by $(\la_i,x_r)$. Then we can assume that in the path representations for both representations the 3 by 3 blocks of $S_3$, as well as the diagonal entries of the 2 by 2 blocks coincide, by Theorem \ref{AB2theorem}. One can now use the braid relations to show that the remaining entries for $S_3$ are determined by this. For example, one can calculate the $((\la_3,x_1),(\la_3,x_1))$ entry of $S_3$ by comparing the $((\la_1,x_3),(\la_3,x_3))$ entries of $S_2S_3S_2$ and $S_3S_2S_3$. The other entries of $S_3$ can be calculated similarly. 

One can similarly prove the uniqueness of the 5- and 4-dimensional simple representations. The four 3-dimensional representations $\{\la_1\la_2\la_3\}$ and
$\{\la_i|\la_j\la_k\}$ are isomorphic as $\K_3$-modules. But their path bases are also  bases of weight vectors, which determines the diagonal matrix $S_3$.

\subsection{Exceptional $B_4$ representations}\label{sec:exc}  We can now list  all new simple representations of $\K_4$ if the eigenvalues of $\sigma_i$ satisfy exactly one irreducible polynomial in Table 2.  Here $\p$ is the prime ideal generated by one of those polynomials, and $Q(R/\p)$ is the corresponding quotient field.
\[
\begin{tabular}{|c|c|c|c|c|c|} $\p$& dim & notation& $\Delta_4^2$& weights& $\K_3$ module \\
\hline
&& & & &\\
$(\la_1\pm i\la_2)$& 2 & $\{\la_1\la_2\}^*$&$ -\la_1^6\la_2^6
$& $(1,2)$& $\{\la_1\la_2\}$\\
$(\la_1+\la_3)$&3&$\{\la_2|\la_1\la_3\}$&
$\la_1^4\la_2^4\la_3^4$&  (2,2), (1,3), & $\{\la_1\la_2\la_3\}$ \\
$(\la_2+\la_3)$&4&$\{\la_1^2\la_2\la_3\}$&
$-\la_1^6\la_2^3\la_3^3$& (1,2), (1,3)&$\{\la_1\la_2\la_3\}+\{\la_1\}$\\
$(\la_1^3-\la_2^2\la_3)$& 5 & $\{\la_2^2|\la_1^2\la_3\}$ & $\la_1^6\la_2^4\la_3^2$& (2,2), (1,2), (1,3)&$\{\la_1\la_2\la_3\}+\{\la_1\la_2\}$\\
$ (\la_1^2-\theta\la_2\la_3)$&7 & $\{\la_1^3\la_2^2\la_3^2\}$ & $\theta\la_1^4\la_2^4\la_3^4$ & $(1,1),(1,2),(1,3), (2,3)$& see \ref{7dimred}\\
\end{tabular}
\]
\begin{equation}\label{7dimred}
\{\la_1^3\la_2^2\la_3^2\}\ \cong\ \{\la_1\la_2\la_3\}+\{\la_1\la_2\} + \{\la_1\la_3\}.
\end{equation}

\centerline {Table 3}
 As in Table 1, we refer to the representation by the determinant of one  of the standard generators $\sigma_i$, where we have additional symbols if the same determinant has appeared before for a generic representation.
We do not list quotients that can be obtained by permuting
the eigenvalues, and we only list one of the two weights $(i,j)$ and $(j,i)$ for $i\neq j$.

\subsection{Blocks}\label{blocks:sec} We now  give an explicit description of the blocks of $\K_4$ if the eigenvalues $\la_i$ are mutually distinct.
We will list the corresponding exact sequences in the next section. This will be proved in Theorem \ref{oneideal}.

(1) $\la_i=-\la_j$: $[\{\la_i^2\la_k\}, \{\la_i^3\la_k^2\la_j\}, \{\la_j^3\la_k^2\la_i\}, \{\la_j^2\la_k\}]$, 
$[ \{\la_i\la_k\}, \{\la_k^3\la_i^2\la_j\}, \{\la_k^3\la_j^2\la_i\}, \{\la_j\la_k\}]$,


(2) $\la_i=-\theta\la_j$:  $[ \{\la_i\}, \{\la_i\la_j\}, \{ \la_j\}]$,  $[\{\la_i^2\la_k\}, \{\la_1^3\la_2^3\la_3^3\}_\theta, \{\la_j^3\la_k^2\la_i\}]$ and

\hskip 6.8em   $[\{\la_j^2\la_k\}, \{\la_1^3\la_2^3\la_3^3\}_{\theta^2}, \{\la_i^3\la_k^2\la_j\}]$,

(3) $\la_i^2=-\la_j^2$: $[ \{\la_i^2\la_j\}, \{\la_j^2\la_i\}, \{\la_i\}, \{\la_j\}]$

(4) $\la_i^2=-\la_j\la_k$: $[\{\la_i\}, \{\la_i\la_j\la_k\}, \{\la_j\la_k\}]$, $[\{\la_j^2\la_k\}, \{\la_j^3\la_i^2\la_k\}, \{\la_i^2\la_j\}]$,


(5) $\la_i^2=\theta\la_j\la_k$: $[\{\la_j\la_k\}, \{\la_1^3\la_2^3\la_3^3\}_\theta, \{\la_i^4\la_j^2\la_k^2\}, \{\la_i\}]$,

(6) $\la_i^3=\la_j^2\la_k$: $[\{\la_i\},  \{\la_i^3\la_j^2\la_k\}, \{\la_j^4\la_i^2\la_k^2\}, \{\la_j^2\la_k\}]$

\subsection{Exact sequences}\label{exact:sec} We list exact sequences corresponding to the blocks in the previous section.  Here the notation is somewhat ambiguous
in case the listed representation is not simple: we may have to replace it by its transpose representation, see
Remark \ref{pathunique}. For the explicit regular representations referred to before, either the given sequence or the sequence with the arrows
reversed for the transposed representations will hold. The first blocks in (2) and (4) appear for representations where $\sigma_1$ and $\sigma_3$ have the same image; the corresponding exact sequences have already appeared for $\K_3$, see Section \ref{K3:sec}, and are not listed here again. Also, we have only listed one exact sequence for the second and third block listed in (2), as they are very similar.
\vskip .2cm
\[
\begin{tabular}{|c|c|c|c|c|c|} $\p$& exact sequence for $Q(R/\p)$ \\
\hline&\\
$(\la_1+\la_3)$&$0 \to\{\la_1^2\la_2\}\to \{ \la_1^3\la_2^2\la_3\}\to\{\la_1\la_2^2\la_3^3\}\to\{\la_2\la_3^2\}\to 0$\\
$(\la_2+\la_3)$&$0 \to\{\la_1\la_2\}\to \{ \la_1^3\la_2^2\la_3\}\to\{\la_1^3\la_2\la_3^2\}\to\{\la_1\la_3\}\to 0$\\
$(\la_1+\theta\la_2)$&$0\to  \{\la_1^2\la_3\}\to\{\la_1^3\la_2^3\la_3^3\}_\theta\to\{\la_1\la_2^3\la_3^2\}\to 0$\\
$(\la_1\pm i\la_2)$&$0\to \{\la_1\}\to \{\la_1^2\la_2\}\to\{\la_1\la_2^2\}\to\{\la_2\}\to 0$\\
$(\la_1^2+\la_2\la_3)$&$ 0\to \{\la_2^2\la_3\}\to \{\la_1^2\la_2^3\la_3\}\to \{ \la_1^2\la_2\}\to 0$\\
$ (\la_1^2-\theta\la_2\la_3)$&$ 0\to \{\la_2\la_3\}\to \{\la_1^3\la_2^3\la_3^3\}_\theta\to \{ \la_1^4\la_2^2\la_3^2\}\to\{\la_1\}\to 0$\\
$(\la_1^3-\la_2^2\la_3)$&$0\to \{\la_2^2\la_3\}\to \{\la_1^2\la_2^4\la_3^2\}\to \{\la_1^3\la_2^2\la_3\}\to \{\la_1\}\to 0$ \\
\end{tabular}
\]
\vskip .2cm
The following theorem gives proofs and concise statements of the results listed in the previous subsections.

\begin{theorem}\label{oneideal}
Let $\p$ be a prime ideal generated by one of the polynomials listed in Table 2, and let $F=Q(R/\p)$ be the quotient field for $R/\p$. Then we have

(a) The blocks of $\K_4$, defined over $F$, are given by equivalence classes of regular modules, where $V_\ga\sim V_{\ga'}$ if and only if $\la_{\ga}=\la_{\ga'}$ mod $\p$ for the eigenvalues $\la_{\ga}, \la_{\ga'}$ of $\Delta_4^2$. See Section \ref{blocks:sec} for explicit lists.

(b) Order the labels $\ga$ in a block by alphabetical order for the eigenvalues $\la_\ga$. After possibly replacing $V_\ga$ by $V_\ga^t$ if necessary, we obtain exact sequences
$$0\ \to\  V_{\ga_1}\ \to\  V_{\ga_2}\ \to\  ...\  \to\  V_{\ga_r}\ \to\  0,$$
where the induced subquotient modules for each $V_{\ga_i}$ are simple. See the table above for explicit exact sequences.
\end{theorem}

$Proof.$ We have determined composition series with simple factors for all regular modules in the previous lemmas. Isomorphisms between the various subquotient modules have been established in Lemma \ref{uniqueness}. Replacing a module $V$ by $V^t$, we can interchange sub and quotient modules, if necessary. So we obtain the claimed exact sequences.

\subsection{No more exceptional representations}\label{nomore:sec} We have listed new simple representations of $\K_4$ when the eigenvalues satisfy
one of the relations in   Proposition \ref{semisimplevalues:prop}. These representations will in general not stay simple if they also satisfy a second equation.
Nevertheless, we have the following result.

\begin{theorem}\label{nomoredecomposition} Table 1 and Table 2 list all possible
simple representations of $\K_4$ for which $\la_i\neq \la_j$, $1\leq i\neq j\leq 3$. The simple modules in a composition series
of a regular module in case when the eigenvalues satisfy two of the relations in Proposition \ref{semisimplevalues:prop}  are listed in Table 4.
\end{theorem}
\vskip .2cm
\[
\begin{tabular}{|c|c|c|} module&  ideal&  simple components\\
\hline&&\\
$\{\la_1^3\la_2^2\la_3\}$&$(\la_1+\la_3, \la_1^3-\la_2^2\la_3$ )&$\{\la_2|\la_1\la_3\}+\{\la_1\la_2\}^*+\{\la_1\}$\\
$\{\la_1^3\la_2^2\la_3\}$&$(\la_2+\la_3, \la_1^3-\la_2^2\la_3)$&$\{\la_1^2\la_2\la_3\}+\{\la_1\}+\{\la_2\}$\\
$\{\la_1^3\la_2^2\la_3\}$&$(\la_2^2+\la_1\la_3, \la_1^3-\la_2^2\la_3)$&$\{\la_2^2\la_1\}+\{\la_1\la_3\}^*+\{\la_1\}$\\
$\{\la_1^4\la_2^2\la_3^2\}$&$(\la_2^3-\la_1^2\la_3, \la_3^3-\la_1^2\la_2,\la_2^2+\la_3^2)$&$\{\la_1^2\la_2\}+\{\la_1^2\la_3\}+\{\la_2\la_3\}^*$\\
$\{\la_1^4\la_2^2\la_3^2\}$&$(\la_2^3-\la_1^2\la_3, \la_3^3-\la_1^2\la_2,\la_2+\la_3)$&$\{\la_1|\la_2\la_3\}+\{\la_1\la_2\}^*+\{\la_1\la_3\}^*+\{\la_1\}$\\
$\{\la_1^4\la_2^2\la_3^2\}$&$(\la_2^3-\la_1^2\la_3, \la_1^2-\theta\la_2\la_3, \la_1^2+\la_3^2)$&$ \{\la_1^2|\la_2^2\la_3\}+\{\la_1\la_3\}^*+\{\la_1\}$\\
$\{\la_1^4\la_2^2\la_3^2\}$&$(\la_2^3-\la_1^2\la_3, \la_1^2-\theta\la_2\la_3,\la_1+\la_3)$&$ \{\la_2^2\la_1\la_3\}+\{\la_1^2\la_3\}+\{\la_1\}$\\
$\{\la_1^3\la_2^3\la_3^3\}_\theta$&$ (\la_1+\theta\la_2, \la_1^2-\theta\la_2\la_3)$&$\{\la_1^2\la_3\}+\{\la_2\la_3\}+\{\la_2^2\la_1\la_3\}$\\
$\{\la_1^3\la_2^3\la_3^3\}_\theta$&$ (\la_2+\theta\la_3, \la_1^2-\theta\la_2\la_3)$&$\{\la_1^2|\la_3^2\la_2\}+\{\la_1\la_2\}^*+\{\la_2\}+\{\la_3\}$ \\
$\{\la_1^3\la_2^3\la_3^3\}_\theta$&$ (\la_1+\theta\la_2, \la_2+\theta\la_3)$&$\{\la_1^2\la_3\}+\{\la_2^2\la_1\}+\{\la_3^2\la_2\}$ \\
\end{tabular}
\]
\vskip .2cm
\centerline{Table 4}
\vskip .2cm
$Proof.$ Let $V$ be a regular module. Then we have a canonical basis of weight vectors. 
 It follows from Theorem \ref{semisimple:thm} that $V$ is simple over the quotient field
$Q(R/\q)$ for any prime ideal $\q$ of $R$ which does not contain an ideal listed in Table 2.
If it only contains one of these ideals, we have a composition series whose factors can be determined in the table in Section \ref{exact:sec}.
We now apply the same strategy to these modules for prime ideals $\q$ which contain more than one of the ideals.
Picking two such ideals $\p_1$ and $\p_2$, we observe the following scenarios:
\vskip .2cm
(a) There is no proper ideal which contains both $\p_1$ and $\p_2$.

(b) If $\p_1$ and $\p_2$ generate a proper ideal, determine all maximal ideals $\m$  which contain both $\p_1$ and $\p_2$. For determining 
all possible subquotient modules, we can already use the information obtained for both $R/\p_1$ and $R/\p_2$, see examples below.
\vskip .2cm
Case (a) can happen when a polynomial $\la_i-\la_j$ is in the ideal generated by $\p_1$ and $\p_2$; recall that we assume
these polynomials to be invertible in $R$. These two ideals likely would still be contained in a nontrivial ideal without our assumption
$\la_i\neq \la_j$. 
 As an example, consider the two ideals $\p_1=(\la_1^2-\theta\la_2\la_3)$ and $\p_2=(\la_2^2-\theta\la_1\la_3)$ for $V=\{\la_1^3\la_2^3\la_3^3\}_\theta$. Then we obtain
$$\la_1^3-\la_2^3=\la_1(\la_1^2-\theta\la_2\la_3)-\la_2(\la_2^2-\theta\la_1\la_3)\ \in \ (\p_1,\p_2).$$
As  $(\la_1-\la_2)$ is invertible, one of the two factors of $(\la_1-\theta\la_2)(\la_1-\theta^2\la_2)=(\la_1^3-\la_2^3)/(\la_1-\la_2)$ has to be in in the prime ideal  $(\p_1,\p_2)$. If $(\la_1-\theta\la_2)\in (\p_1,\p_2)$, then so is $\la_1-\la_3=\la_2^{-1}\theta^{-1}(\la_1^2-\theta\la_2\la_3-\la_1(\la_1-\theta\la_2))$.
One similarly shows that $\la_2-\la_3\in (\p_1,\p_2)$ if $(\la_1-\theta^2\la_2)\in (\p_1,\p_2)$.

For $V=\{\la_1^3\la_2^2\la_3\}$, the ideal generated by  $(\la_1+\la_3)$ and $(\la_2+\la_3)$ obviously contains the invertible polynomial $\la_1-\la_2$. 
As another example for case (a), one similarly shows that the ideal
generated by $(\la_2^2+\la_1\la_3)$ and $(\la_2+\la_3)$ contains $\la_1-\la_2$ from
$$\la_1-\la_2=\la_3^{-1}[(\la_2^2+\la_1\la_3)-\la_2(\la_2+\la_3)].$$
For case (b), observe that both $\p_1$ and $\p_2$ define a partition of the weight basis.
If all weights have multiplicity 1, the intersections of these subsets give us a refined partition. It only remains to check that
the path vectors labeled by the elements of any of its components do indeed span a simple subquotient of $V$.
As an example for this case, consider the ideal generated by $(\la_1^2-\theta\la_2\la_3)$ and $(\la_2+\theta\la_3)$
for the 9-dimensional generic representation $\{\la_1^3\la_2^3\la_3^3\}_\theta$.
Then it follows from Tables 1 and 3 that the equation $\la_2=-\theta\la_3$
induces  a partition
of its weights into $\{ (2,2), (1,2), (2,1)\}$ 
and the remaining weights. Similarly, the equation $\la_1^2=\theta\la_2\la_3$ induces a partition of the weights of 
 $\{\la_1^3\la_2^3\la_3^3\}_\theta$ into $\{ (2,2), (3,3)\}$ and its complement.
  Hence, if both
equations are satisfied, we obtain the refined partition
$$\{ (2,2)\}\ \cup\ \{ (1,2), (2,1)\}\ \cup\ \{ (3,3)\}\ \cup\ \{(1,1),  (1,3), (3,1), (2,3), (3,2)\}.$$
We deduce from this and Table 3 the result stated in the second to last line of Table 4.
Observe that the corresponding simple subrepresentations have already appeared before. E.g. one can check directly that 
our two defining equations also imply $\la_1^2=-\la_2^2$, a necessary condition for the $\K_4$ representation $\{\la_1\la_2\}^*$.
We similarly obtain the results listed in Table 4 for the other cases concerning the 6- and 8-dimensional regular modules.

The situation is slightly more complicated for the 8-dimensional representation, say, $\{\la_1^4\la_2^2\la_3^2\}$. As the weight (1,1) has multiplicity 2,
there is more than one refinement of our path basis compatible with the partitions coming from the ideals $\p_1$ and $\p_2$. They can be associated to different
ideals $\q_1$ and $\q_2$, both of which contain $\p_1$ and $\p_2$.
As an example, consider the two ideals
$\p_1=(\la_2^3-\la_1^2\la_3)$ and $\p_2=(\la_3^3-\la_1^2\la_2)$.
First, observe that  we have subquotients isomorphic to $\{\la_1^2\la_2\}$ and $\{\la_1^2\la_3\}$, by Lemma \ref{8dimcomplem}. Moreover, we obtain
$$(\la_2+\la_3)(\la_2^2+\la_3^2)=(\la_2-\la_3)^{-1}[\la_2(\la_2^3-\la_1^2\la_3)-\la_3(\la_3^3-\la_1^2\la_2)]\ \in\  (\p_1, \p_2).$$
We consider the two cases of ideals  $\q_1=(\la_2^3-\la_1^2\la_3, \la_3^3-\la_1^2\la_2, \la_2+\la_3)$ and $\q_2=(\la_2^3-\la_1^2\la_3, \la_3^3-\la_1^2\la_2, \la_2^2+\la_3^2)$.
They are both nontrivial, with $(i,1,-1)$ and $(\sqrt{i},1,-i)$ elements in their loci, respectively. 

In the first case, we derive from $\la_2=-\la_3$ and the relations we started out with that $\la_1^2=-\la_2^2=-\la_3^2$.
 This splits the subquotients $\{\la_1^2\la_i\}=\{\la_1\la_i\}^*+\{\la_1\}$ 
for $i=2,3$. The remaining weights are (23), (32) and, possibly (1,1) if the $\{\la_1\}$ summands in the split of $\{\la_1^2\la_i\}$ coincide for $i=2,3$. As $\la_2^2\neq -\la_3^2$, there is no 2-dimensional representation with weights (2,3) and (3,2).
Hence, the remaining summand is $\{\la_1|\la_2\la_3\}$,  which gives us line 5 in Table 4 for the representation $\{\la_1^4\la_2^2\la_3^2\}$.

In the second case, we deduce from $\la_2^2=-\la_3^2$ and $\la_2^3=\la_1^2\la_3$ that $\la_1^2=-\la_2\la_3\neq -\la_3^2$,
and similarly $\la_1^2\neq -\la_2^2$.
Hence, the subquotients $\{\la_1^2\la_2\}$ and $\{\la_1^2\la_3\}$ remain irreducible, 
and the remaining summand necessarily has to be $\{\la_2\la_3\}^*$. This  gives us line 4 in Table 4 for the representation $\{\la_1^4\la_2^2\la_3^2\}$.

If $\p_1=(\la_2^3-\la_1^2\la_3)$ and $\p_2=(\la_1^2-\theta\la_2\la_3)$, we similarly determine ideals $\q_1=(\p_1,\p_2,\la_1^2+\la_3^2)$ and  $\q_2=(\p_1,\p_2,\la_1+\la_3)$,
which lead to subquotients listed 
in lines 6 and 7 of Table 4.
\section{Related results and applications} 

\subsection{Schur elements and semisimplicity}\label{Schursec}  Our Theorem A proves a special case of a general conjecture, which would determine the exact values of parameters for which a cyclotomic Hecke algebra is semisimple. We sketch this here, following the review in \cite{CJ}, Chapter 2. 
A symmetrizing trace for an algebra $\A$ over a ground ring $R$  is a functional $tr:\A\to R$ such that
$(a,b)=tr(ab)$ defines a symmetric non-degenerate bilinear form. A canonical symmetrizing trace for a cyclotomic Hecke algebra has been conjectured by Brou\' e, Malle and Michel for all cyclotomic Hecke algebras.
This has been verified for a large class of such Hecke algebras, see the review in \cite{CJ}, but not for $G_{25}$. 
The trace can be characterized by so-called Schur elements, which have been calculated for all cyclotomic Hecke algebras in \cite{Malledegree},
assuming the existence of the trace. It is well known that an algebra $\A$ is semisimple if and only if all Schur elements
are nonzero. Comparing the Schur elements for $G_{25}$ in \cite{Malledegree} with the polynomials in Theorem A shows that we have confirmed the conjecture for this particular Hecke algebra.

\subsection{Reconstructing categories of type $G_2$} Our original motivation for this research comes from our project of classifying semisimple braided rigid tensor categories $\Ca$ whose tensor product rules are the ones of the representations of the Lie algebra of type $G_2$ or its associated fusion tensor categories. Such classifications have been obtained for classical Lie types in \cite{KW} and \cite{TW2}, but are not known for exceptional Lie types. In \cite{MW2} it was shown that if
the image of the braid group $B_3$ in $\End(V^{\otimes 3})$ generates this algebra, then also $\End(V^{\otimes n})$ is generated by the image of $B_n$; here $V$ is the object in $\Ca$ corresponding to the simple 7-dimensional representation of the Lie algebra $\g(G_2)$.  From this we can conclude that the given tensor category $\Ca$ must be equivalent to Rep $U_q\g(G_2)$, the Drinfeld-Jimbo qantum group of type $G_2$, with $q$ not a root of unity.

In order to prove the missing step, we need to examine all possible representations of the braid group $B_4$ into $\End(V^{\otimes 4})$ to eventually prove that the image of $B_3$ must generate $\End(V^{\otimes 3})$.
This will be done in \cite{MW3}, also for the associated finite fusion categories of type $G_2$, using the results of this paper. It is expected that the results of our paper would also be useful in proving similar classification results for braided tensor categories of other exceptional Lie types.

\subsection{Scope of our method} A priori, our method of using path bases would only work if both $\K_2$ and $\K_3$ are semisimple.  We managed to make it work even in the case where $\K_3$ is not semisimple. It might be possible to extend this even to the general case with also $\K_2$ not semisimple, i.e. the braid generators are not diagonalizable. Similar situations were studied in \cite{GW} in the context of Iwahori-Hecke algebras of type $A$. But more work would be needed.

We expect that our methods will also be useful for determining the structure of other cyclotomic Hecke algebras corresponding to complex reflection groups where some of the generators have order 3, such as e.g. $G_{32}$.

\end{document}